\documentclass[a4paper,reqno,11pt]{amsart}
\usepackage{mathrsfs}
\usepackage[all]{xy}
\usepackage{txfonts}
\usepackage{amsfonts}
\usepackage{amssymb}
\usepackage{amsfonts,amsmath,amsthm,amssymb,stmaryrd}
\usepackage{hyperref}
\usepackage[numbers,sort&compress]{natbib}

\newtheorem{Theorem}{Theorem}[section]
\newtheorem{Lemma}{Lemma}[section]
\newtheorem{Proposition}{Proposition}[section]

\numberwithin{equation}{section}
\usepackage[left=1 in, right=1 in,top=1 in, bottom=1 in]{geometry}

\def\XXint#1#2#3{{\setbox0=\hbox{$#1{#2#3}{\int}$ }
\vcenter{\hbox{$#2#3$ }}\kern-.6\wd0}}

\DeclareMathOperator{\diver}{div}

\def\r3{\mathbb{R}^3}

\newcommand{\pa}{\partial}

\makeatletter
\@namedef{subjclassname@2020}{\textup{2020} Mathematics Subject Classification}
\makeatother

\allowdisplaybreaks
\begin{document}
\bibliographystyle{plain}

\title[D\MakeLowercase{ecay estimates for the compressible viscoelastic equations}]{D\MakeLowercase{ecay estimates for the compressible viscoelastic equations in an exterior domain}}

\author{Jieling Deng}
\address{School of Mathematical Sciences, South China Normal University, Guangzhou, Guangdong 510631, China.}
\email[J. Deng]{jieling202101@163.com}

\author{Yong Wang}
\address{School of Mathematical Sciences, South China Normal University, Guangzhou, Guangdong 510631, China.}
\email[Y. Wang]{wangyongxmu@163.com}

\author{Jianquan Yang}
\address{School of Mathematical Sciences, South China Normal University, Guangzhou, Guangdong 510631, China.}
\email[J. Yang]{yangjianquanscnu@163.com}

\thanks{Corresponding author: Yong Wang (wangyongxmu@163.com).}

\begin{abstract}
In this paper, we study the compressible viscoelastic equations in an exterior domain. We prove the $L_2$ estimates for the solution to the linearized problem and show the decay estimates for the solution to the nonlinear problem. In particular, we obtain the optimal decay rates of the solution itself and its spatial-time derivatives in the $L_2$-norm.
\end{abstract}
\keywords{Compressible viscoelastic flows; Strong solutions; Exterior problems; Decay rates.}
\subjclass[2020]{76A10; 35Q35; 35G55.}

\maketitle

\section{Introduction}
In this paper, we consider the initial-boundary value problem (IBVP) to the three-dimensional compressible viscoelastic equations: for $( t,x ) \in \mathbb{R}^+ \times \Omega$,
\begin{align}\label{1.1'}
\begin{cases}
	\rho _t+\mathrm{div}( \rho u ) =0,\\
	( \rho u ) _t+\mathrm{div}( \rho u\otimes u ) -\mu \Delta u-( \lambda +\mu ) \nabla \mathrm{div}u+\nabla P( \rho ) =\mathrm{div}( \rho FF^T ),\\
	F_t+u\cdot \nabla F=\nabla uF,\\
\left.	( \rho ,u,F )\right| _{t=0}=( \rho _0,u_0,F_0 ) ,\\
\end{cases}
\end{align}
which is subject to the boundary condition
\begin{align}\label{bound}
u|_{\partial \Omega}=0,\quad t>0.
\end{align}
Here $\Omega\subset\mathbb{R}^3$ is the exterior of a bounded domain with smooth boundary $\pa\Omega$. Here $\rho \in  \mathbb{R} $, $u\in \mathbb{R} ^3$, $F\in M^{3\times 3}$(the set of $3\times 3$ matrices with positive determinants) denote the density, th velocity, and the deformation gradient, respectively. $\mu , \lambda $ are the shear viscosity and the bulk viscosity coeffcients of the fluid, respectively, which are assumed to satisfy the following physcial condition:
$$\mu >0,\quad 2\mu +3\lambda \geqq 0.$$
The pressure term $P(\rho)$ is a suitably smooth function of $\rho$ for $\rho>0$. The symbol $\otimes $ denotes the Kronecker tensor product, $F^T$ means the transpose matrix of $F$, and the notation $u\cdot \nabla F$ is understood to be $(u\cdot \nabla) F$. The positive parameter $\alpha$ represents the speed of propagation of shear waves. \par
System \eqref{1.1'} is one of the basic macroscopic model for compressible flows exhibiting elastic properties, which corresponds to the so-called Hookean linear elasticity. For more physical background about this system, we refer to \cite{Gurtin1981,Lin-Liu-Zhang2005}\par
In the context of hydrodnamics,  the motion of the fluid is classically described by a time-dependent family of orientation preserving diff-isomorphisms $x(t,X),0\leqq t<T$. Material points X in the reference configuration are deformed to the spatial position $x(t,X)$ at the time $t$. Then the deformation tensor $\widetilde{F}(t,X)$ is defined as
$$\widetilde{F}(t,X)\overset{\mathrm{def}}{\operatorname*{\operatorname*{=}}}\frac{\partial x(t,X)}{\partial X}.$$
When working in the Eulerian coordinates, we defined $\widetilde{F}(t,X)$ is defined as
$$F\left(t,x(t,X)\right)=\widetilde{F}(t,X).$$
Applying the chain rule, we see that $F(t,X)$  satisfies the following transport equation
$$F_t+u\cdot\nabla F=\nabla uF,$$
which stands for
$$F_{t}^{ij}+u^{k}\nabla_{k}F^{ij}=\nabla_{k}u^{i}F^{kj},\quad
\mathrm{for}\ i,j=1,\ldots,N,$$
where
$$\nabla_i=\frac{\partial}{\partial x_i},\quad F^{ij}=\frac{\partial x^i}{\partial X^j},\quad(\nabla u)_{ij}=\nabla_ju^i.$$
We take the convention of summing over repeated indices. We also assume that
\begin{align}\label{1.2.}
\mathrm{div}(\rho F^T)=0,\quad F^{lk}(0)\nabla_{l}F^{ij}(0)=F^{lj}(0)\nabla_{l}F^{ik}(0).
\end{align}
It is standard that the condition \eqref{1.2.} is preserved by the flow, which has been proved in \cite{Hu-Wang2011}.

There are lots of literature on the well-posedness of solutions of the Cauchy or initial-boundary value problem of the compressible viscoelastic system \eqref{1.1'}.
For the 3D Cauchy problem, Hu and Wang \cite{Hu-Wang2010} established the existence and uniqueness of the local strong solution with large initial data in Sobolev spaces.
Under the condition that the initial data were close to the constant equilibrium state, Hu and Wu \cite{Hu-Wu2013} showed the global existence and optimal decay rates of strong solutions in the $H^2$ framework.
For the multi-dimensional case, Qian-Zhang \cite{Qian-Zhang2010} and Hu-Wang \cite{Hu-Wang2011} obtained the global-in-time existence of the strong solution with initial data close to an equilibrium state in Besov spaces.
Later, Zhu \cite{Zhu2022} proved the global small solutions of the 3D compressible viscoelastic system without using any structural assumption, which played a crucial role in previous analysis mentioned above; see also \cite{Pan-Xu-Zhu2022,Liu-Meng-Wu-Zhang2023} for the framework of critical Besov spaces.
For the initial-boundary value problem, the global-in-time existence of the strong solution near equilibrium was established by Qian in \cite{Qian2011}, where the cases of half space, bounded domain and exterior domain were considered; see also \cite{Hu-Wang2015,Chen-Wu2018}.
For more research work on the well-posedness theory of compressible viscoelastic system, one can refer to \cite{ref2020.1,Cui-Hu2022,Han-Zi2020,Tan-Wang-Wu2020,Wu-Wang2023,Wang-Xie2023,Gu-Wang-Xie2024,Huang-Yao-You2025} and the references cited therein.

In addition to well-posedness, the large time behavior of the solutions is another important subject to study.
For the three-dimensional Cauchy problem, by using the Fourier splitting method and the Hodge decomposition, Hu-Wu \cite{Hu-Wu2013} showed the following optimal decay estimates of the strong solution: for $2\leqq p\leqq6$,
    \begin{equation*}
        \begin{aligned}
        &\|(\rho-1,u,F-I)(t)\|_{L_p}\leqq C(1+t)^{-\frac{3}{2}(1-\frac{1}{p})}, \\
        &\|\nabla(\rho-1,u,F-I)(t)\|_{H^1}\leqq C(1+t)^{-\frac{5}{4}},
        \end{aligned}
    \end{equation*}
provided that the initial perturbation is sufficiently small in $H^2$ and belongs to $L^1$.
Later, Li-Wei-Yao \cite{Li-Wei-Yao2016} extended the result in \cite{Hu-Wu2013} by deriving the following optimal time decay estimates on the higher-order derivatives of the classical solution: for $2\leqq p\leqq\infty$ and $k=0,1,\ldots,l-1$,
    \begin{equation*}
        \begin{aligned}
        &\|\nabla^k(\rho-1,u,F-I)(t)\|_{L_p}\leqq C(1+t)^{-\frac{k}{2}-\frac{3}{2}(1-\frac{1}{p})},\\
        &\|\nabla^k(\rho-1,u,F-I)(t)\|_{H^{l-k}}\leqq C(1+t)^{-\frac{3}{4}-\frac{k}{2}},
        \end{aligned}
    \end{equation*}
provided that the initial perturbation is sufficiently small in $H^l$ $(l\geqq3)$ and belongs to $L^1$.
By using a pure energy method, Wu-Gao-Tan \cite{Wu-Gao-Tan2017} also obtained the optimal decay rates of higher-order derivatives of the classical solution: for $k=0,1,\ldots,l-1$,
    \begin{equation*}
        \|\nabla^k(\rho-1,u,F-I)(t)\|_{H^{l-k}}\leqq C(1+t)^{-\frac{k+s}{2}},
    \end{equation*}
provided that the initial perturbation is sufficiently small in $H^l$ $(l\geqq3)$ and belongs to $\dot{H}^{-s}$ or $\dot{B}_{2,\infty}^{-s}$ $(s\in(0,\frac{3}{2}])$.
This relaxed the $L^1$ initial assumption proposed in \cite{Hu-Wu2013,Li-Wei-Yao2016} since $L^1\subset \dot{B}_{2,\infty}^{-\frac{3}{2}}$.
In \cite{Fu-Huang-Wang2023}, Fu-Huang-Wang proved the optimal time decay rate of the highest-order derivative of the classical solution: for $k=0,1,\ldots,l$,
    \begin{equation*}
        \|\nabla^k(\rho-1,u,F-I)(t)\|_{H^{l-k}}\leqq C(1+t)^{-\frac{k+s}{2}},
    \end{equation*}
provided that the initial perturbation is sufficiently small in belonged to $H^l$ $(l\geqq3)$ and belongs to $\dot{H}^{-s}$ $(s\in[0,\frac{5}{2}))$.
It should be noted that Ishigaki \cite{Ishigaki2020} clarified the diffusion wave phenomena caused by the sound wave and the elastic shear wave, and showed that the global classical solution satisfied the following $L_p$ decay estimates: for $1<p\leqq\infty$,
    \begin{equation*}
        \|(\rho-1,u,F-I)(t)\|_{L_p}
        \leqq C(1+t)^{-\frac{3}{2}(1-\frac{1}{p})-\frac{1}{2}(1-\frac{2}{p})},
    \end{equation*}
provided that the initial perturbation is sufficiently small $H^3$ and belongs to $L^1$.
This improved the decay rates of the $L_p$-norm of the solution obtained in \cite{Hu-Wu2013,Li-Wei-Yao2016} for $p>2$.
One can refer to \cite{Bai-Zhang2023} for the pointwise estimates of the solution for system \eqref{1.1'}.
Recently, Fu-Huang-Jiang \cite{Fu-Huang-Jiang2024} studied the temporal decay rates of the solution to the Cauchy problem for the 3D compressible viscoelastic equations with relatively large elasticity coefficient.
By using the energy method with temporal weights and the spectral analysis, they proved that the $k^{th}$ ($k=0,1$) order spatial derivatives of both the density and deformation perturbations converge to zero in $L_2$-norm at a rate of $(1+t)^{-\frac{3}{4}-\frac{k+1}{2}}$, which improved the result in \cite{Hu-Wu2013}. We also refer the readers to \cite{Jia-Peng-Mei2014,Pan-Xu2019,Jia-Peng2017,Xu-Zhang-Wu-Liu2016} for the large time behavior of the solutions in the critical Besov spaces.
As for the large time behavior of the solutions to the initial-boundary value problem, Chen and Wu \cite{Chen-Wu2018} showed the estimation of the exponential convergence rates of the strong solution; see also \cite{Wang-Shen-Wu-Zhang2022} for the study on compressible viscoelastic fluids with the electrostatic effect.
However, to the best of our knowledge, the large time behavior of the solutions for the compressible viscoelastic system \eqref{1.1'} in an exterior domain has been rarely studied.
In fact, this is the aim of this paper.

\bigskip

\textbf{Notations}:
Before stating our main results precisely, at this point we shall explain our notation. Let $D$ be any domain in $\mathbb{R}^3$. $L_p(D)$ denotes the usual $L_p$ space on $D$ with the norm $\|\cdot\|_{p,D}.$ Put
\begin{align*}
W_p^m(D)=\{u\in L_p(D)\mid\|u\|_{m,p,D}<\infty\},\quad\|u\|_{m,p,D}=\sum_{|\alpha|\leqq m}\|\partial_x^\alpha u\|_{p,D},\\
H^m(D)=W_2^m(D),\quad W_p^0(D)=L_p(D),\quad H^0(D)=L_2(D).
\end{align*}
Denote $\mathcal{W}_p^m(D)=\{{\text{the closure of}\  C_{0}^{\infty}(D) \ \text{in} \ W_p^m(D)}\}$, where $C_{0}^{\infty}(D)$ denotes the set of all $C^{\infty}(D)$ functions whose supports are compact in $D$. Next, set
\begin{align*}
\|u\|_{m,p,D}=\sum_{j=1}^3\|u_j\|_{m,p,D},\quad\|\mathbb{U}\|_{m,p,D}=\sum_{j=1}^5\|u_j\|_{m,p,D},\quad \|F\|_{m,p,D}=\sum_{i,j=1}^3\|F_{ij}\|_{m,p,D}.
\end{align*}
Set $\mathbb{PU}=(0,u,0)^{T}$ and $(\mathbb{I}-\mathbb{P})\mathbb{U}=(\rho,0,F)^{T}$ for $\mathbb{U}=(\rho,u,F).$ Set
\begin{align*}
\begin{aligned}\mathbb{W}_{p}^{k,m}(D)&=\{\mathbb{U}\mid(\mathbb{I}-\mathbb{P})\mathbb{U}\in\mathbb{W}_p^k(D),\mathbb{P}\mathbb{U}\in\mathbb{W}_p^m(D)\},\\\|\mathbb{U}\|_{\mathbb{W}_p^{k,m}}(D)&=\|(\mathbb{I}-\mathbb{P})\mathbb{U}\|_{k,p,D}+\|\mathbb{P}\mathbb{U}\|_{m,p,D}.\end{aligned}
\end{align*}
When $D=\Omega$, we omit the index $\Omega$. Namely,
\begin{align*}
\|\cdot\|_{p,\Omega}=\|\cdot\|_p,\quad\|\cdot\|_{k,p,\Omega}=\|\cdot\|_{k,p},\quad\|\cdot\|_{\mathbb{W}_p^{k,m}(\Omega)}=\|\cdot\|_{\mathbb{W}_p^{k,m}}.
\end{align*}
Set
\begin{align*}
&C^\ell([t_1,t_2],B)\\
=&\{u(t)\mid\ell \text{-times continuously differentiable function of} \ t\in[t_1,t_2]
\text{with valued in a Banach space} B\},
\end{align*}
and
\begin{align*}
L_2((t_1,t_2),B)=\{u(t)\mid L_2\text{-function of}\ t\in[t_1,t_2] \text{with valued in}\  B\}.
\end{align*}
Let
\begin{align*}
&N_{1}(0,\infty)^{2}
=\sup_{0\leqq t<\infty}\left(\|\mathbb{U}(t)-\overline{\mathbb{U}}_{0}\|_{3,2}^{2}
+\|\partial_{t}\mathbb{U}(t)\|_{\mathbb{W}_{2}^{2,1}}^{2}\right),\\
&N_2(0,\infty)^2=\sup_{0\leqq t<\infty}\left(\|\mathbb{U}(t)-\overline{\mathbb{U}}_{0}\|_{4,2}^{2}
+\|\partial_{t}\mathbb{U}(t)\|_{\mathbb{W}_{2}^{3,2}}^{2}
+\|\partial_{t}^{2}\mathbb{U}(t)\|_{\mathbb{W}_{2}^{2,0}}^{2}\right).
\end{align*}
Next, let us explain some necessary condition on the initial data $(\rho_0,u_0,F_0)$. To do this, for a while we assume that $(\rho,u,F)\in X^{k}(0,\infty)$ is a solution of IBVP \eqref{1.1'}. Then, $\partial_{t}^{j}(\rho,u,F)|_{t=0},j\geqq1$, are determined successively by the initial data $(\rho_0,u_0,F_0)$ through Eq.\,\eqref{1.1'}. Namely, if we put $(\rho_{j},u_{j},F_{j})=\partial_{t}^{j}(\rho,u,F)|_{t=0}$ and then $(\rho_{j},u_{j},F_{j})$ is generated inductively from $(\rho_0,u_0,F_0)$ by means of Eq.\,\eqref{1.1'}. From the definition of $X^{k}(0,\infty)$, the solution $(\rho,u,F)$ belonging to $X^{k}(0,\infty)$ should satisfy the following condition at $t=0$:
\begin{align*}
\rho_{0}-1\in H^{k+2},\quad\rho_{j}\in H^{k+2-j},\quad(j=1,k),\\
u_{j}\in H^{k+2-2j}\quad(j=0,1,k),\quad u_{j}|_{\partial\Omega}=0\quad(j=0,1),\\
F_{0}-I\in H^{k+2},\quad F_{j}\in H^{k+2-j}\quad(j=1,k).
\end{align*}
This is a requirement on the initial data $(\rho_0,u_0,F_0)$ in order to obtain solutions belonging to $X^{k}(0,\infty)$. If $(\rho_0,u_0,F_0)$ satisfies these conditions above, we say that $(\rho_0,u_0,F_0)$ satisfies the \textbf{$k^{th}$ order compatibility condition and regularity}.
To denote various constants we use the same letter $C$, and $C(A,B,\ldots)$ denotes the constant depending on the quantities, $A,B,\ldots\mathrm{O}(\cdot)$ means the large order.\par
\bigskip
\textbf{Main contributions}: Our main result is Theorem \ref{the1.2}. To prove this theorem, we require decay estimates for the linearized compressible viscoelastic system in an exterior domains, which will be combined with the Duhamel principle to establish Theorem \ref{the1.2}. Since we are considering an exterior domain problem, when analyzing the linearized compressible viscoelastic system , we must consider for both local energy decay and some estimates in $\mathbb{R}^3$.\par
In Section 3, we cite the work of Ishigaki and Kobayashi \cite{Ishigaki-Kobayashi2025} to obtain relevant properties of the linear viscoelastic system in an exterior domains along with local energy decay estimates. The local energy decay estimates of $\mathbb{U}(t)=\mathbb{T}(t)\mathbb{U}_0$ hold for all $\begin{aligned}m=0,1,2,\ldots\end{aligned}$ and $t\geqq1:$
\begin{align*}
\begin{array}{rcl}\|\partial_t^m\mathbb{T}(t)\mathbb{U}_0\|_{\mathbb{W}_2^{1,2}(\Omega_b)}&\leqq C_{b,b_0,m}\end{array}t^{-2-m}\|\mathbb{U}_0\|_{\mathbb{W}_2^{1,0}(\Omega)},
\end{align*}
where $C_{b,b_0,m}$ is a constant depending only on $b,b_0,m$ .

In Section 4, we investigate the decay estimates of the linear viscoelastic equations in $\mathbb{R}^3$. By employing the Helmholtz decomposition
$$v=-\Lambda^{-1}\nabla d-\Lambda^{-1}\mathrm{div}\omega$$
from Hu and Wu \cite{Hu-Wu2013} , we simplify the eigenvalue computation and derive more comprehensive decay estimates than those obtained by Wu and Hu, including temporal and spatial derivative decay estimates, thereby yielding Theorem \ref{th3.1}.\par
With these preparatory results, in Section 5, we utilize cut-off technique to partition the domain.Put $\begin{aligned}\mathbb{U}_m=(n_m,v_m,E_m)=\mathbb{U}-(1-\varphi)\mathbb{U}_c\end{aligned}$ and $\mathbb{U}_\infty(t)=(n_\infty(t),v_\infty(t),E_\infty(t))=\varphi_\infty\mathbb{U}(t).$ By applying both local energy decay and some estimates in $\mathbb{R}^3$ under appropriate conditions, we establish the decay estimates for the linear viscoelastic system in an exterior domains, leading to Theorem \ref{the1.3}.\par
In Section 6, we focus on the key estimates of nonlinear terms. Building upon the results from Section 5, we complete the proof of Theorem \ref{the1.2}.\par

\section{Main result}
\begin{Lemma}
The density $\rho$ and deformation tensor $F$ satisfy
\begin{align}\label{1.6..}
\mathrm{div}(\rho F^T)=0\quad \rho_0\mathrm{det}F_0=1\quad \mathrm{and}\quad F^{lk}\nabla_lF^{ij}-F^{lj}\nabla_lF^{ik}=0,
\end{align}
if the initial data $(\rho,F)|_{t=0}=(\rho_0,F_0)$ satisfy
\begin{align}\label{1.4.}
\mathrm{div}(\rho_0F_0^T)=0\quad\mathrm{and}\quad F_0^{lk}\nabla_lF_0^{ij}-F_0^{lj}\nabla_lF_0^{ik}=0,
\end{align}
respectively.
\end{Lemma}

Define
\begin{align}\label{1.5.}
G\overset{\mathrm{def}}{\operatorname*{\operatorname*{=}}}\frac{\partial X^{-1}(t,x)}{\partial x}=F^{-1},\quad H\overset{\mathrm{def}}{\operatorname*{\operatorname*{=}}}G-I,
\end{align}
then the second identity in \eqref{1.5.} is equivalent to the curl free property of $G$ or $H$. From $\eqref{1.1'}_3$ we can find that $G$ satisfies
\begin{align*}
G_t+u\cdot\nabla G+G\nabla u=0.
\end{align*}
For further, let $\widetilde{X}(t,x)$ be determined by
\begin{align*}
\begin{cases}
\tilde{X}_t=-u \cdot\nabla\tilde{X},\\\tilde{X}|_{t=0}=\tilde{X}_0.
\end{cases}
\end{align*}
One can verify that $\widetilde{X}(t,x)=\widetilde{X}(0,X^{-1}(t,x))$ and $\nabla_j\widetilde{X}^i$ satisfies the asme equation as $G$.

\begin{Proposition}\label{The1.1}
Let $\Omega$ be a domain in $\mathbb{R}^d, d=2,3$ . Assume that $\partial\Omega\in C^3$ and $(\rho_0,u_0,F_0)\in H^2$. We further assume that the initial data satisfies \eqref{1.4.} and $G_0^{ij}=\nabla_j\widetilde{X}_0^i$ for some $\tilde{X}_{0}$ satisfying $\tilde{X_0}-x=0$ on $\partial\Omega$.  Then, 
there exist an numbers $c_0$ such that, if
\begin{align*}
\|\rho_0-1,u,F_0-I\|_2\leqq c_0,
\end{align*}
then problem \eqref{1.1'} has a unique solution $(\rho,u,F)\in\mathcal{C}([0,+\infty],H^2).$ Moreover, the solution decays exponentially to the equilibrium state $(\bar{\rho},0,I)$ in the case of bounded domain and
\begin{align*}
N_{k}(0,\infty)\leqq C\left\|(\rho_{0},u_{0},F_{0})\right\|_{H^{k+2}},\quad k=1,2.
\end{align*}
\end{Proposition}
The main purpose of the paper is to prove the following theorem concerning the rate of convergence of solutions to IBVP \eqref{1.1'}.
\begin{Theorem}\label{the1.2}
Assume that the initial data $(\rho_{0},u_{0},F_{0})$ satisfy the $2^{th}$ order compatibility condition and regularity and $\nabla \rho_{0}+\mathrm{div}{F_0}^T=0$ . Then there exists an $\epsilon>0$ such that if $\|(\rho_{0}-1,u_{0},F_0-I)\|_{4,2}\leqq\epsilon$,  then the solution $(\rho,u,F)$ of the IBVP \eqref{1.1'}-\eqref{bound} has the following asymptotic behaviours as $t\to\infty$:
\begin{align*}
&\|(\rho-1,u,F-I)(t)\|_{2}=O(t^{-\frac{3}{4}});\\
&\|\partial_{x}(\rho,u,F)(t)\|_{\mathbb{W}_{2}^{1,2}}+\|\partial_{t}(\rho,u,F)(t)\|_{1,2}=O(t^{-\frac{5}{4}});\\
&\|(\rho-1,u,F-I)(t)\|_{\infty}=O(t^{-\frac{3}{2}}).
\end{align*}
\end{Theorem}
In order to prove Theorem \ref{the1.2}, we shall use the decay property of solutions to the corresponding linearized problem. Without loss of generality, we assume $P^{\prime}(1)>0$,and denote $\chi_{0}=(P^{\prime}(1))^{-\frac{1}{2}} $. we used  the condition $\diver(\rho F^T)=0$ for all $t\geqq0$. For $\rho>0$, we denote
$$n(t,x)=\rho(\chi_0^2t,\chi_0x)-1,\quad v(t,x)=\chi_0u(\chi_0^2t,\chi_0x),\quad E(t,x)=F(\chi_0^2t,\chi_0x)-I.$$
If we linearize Eq.\,\eqref{1.1'} at the constant state $(1,0,I)$ and we make some linear transformation of the unknown function. Set $E=(E_1,E_2,E_3)^T$, then we have the following initial boundary value problem of the linear operators:
\begin{align}\label{2.8.}
\begin{cases}
n_t+\mathrm{div}v=0,\\
v_t-\mu\Delta v-(\lambda+\mu)\nabla\mathrm{div}v+\nabla n-\mathrm{div}E^{T}=0,\\
E_t-\nabla v=0
\end{cases}
\end{align}
with the boundary and initial condition
\begin{align}\label{2.10..}
v|_{\partial\Omega}=0,\quad(n,v,E)|_{t=0}=(n_0,v_0,E_0)
\end{align}
for $( t,x ) \in \left[ 0,+\infty \right) \times \Omega$. Under the constraints obtained by linearizing \eqref{1.6..} without its second constraint:
$$\nabla n_0+\mathrm{div}{E_0}^T=0,\quad
\partial_{x_l}E_0^{jk}=\partial_{x_k}E_0^{jl}\ (j,k,l=1,2,3).$$
Here $\mathbb{A}$ is the linearized operator given by
\begin{align*}
\mathbb{A}=
\begin{pmatrix}0&\mathrm{-div}&0\\
-\nabla&\mu\Delta+(\lambda+\mu)\nabla\mathrm{div}&\mathrm{div}\\
0&\nabla&0\end{pmatrix}.
\end{align*}
Then, \eqref{2.8.} is wriiten in the form:
\begin{align}\label{2.9.}
\mathbb{U}_t-\mathbb{A}\mathbb{U}=0\quad\mathrm{~for~}t>0,\quad\mathbb{U}|_{t=0}=\mathbb{U}_0,
\end{align}
where $\mathbb{U}_0=(n_0,v_0,E_0)$ and $\mathbb{U}=(n,v,E)$. We consider the
case that there exists a vector field $\psi$ such that $E$ is written as $E=\nabla\psi$ , which automatically satisfies
the identity of \eqref{2.10..}.
In order to prove Theorem \ref{the1.2}, we have to obtain the suitable decay property of solutions to \eqref{2.9.}. In this paper, we shall show the following theorem concerning the decay rate of $\{e^{t\mathbb{A}}\}_{t\geqq0}.$
\begin{Theorem}\label{the1.3}
Put $\mathbb{U}(t)=e^{t\mathbb{A}}\mathbb{U}_0$.
Let $1\leqq q\leqq2$. All the constants $C$ in the theorem depend only on $q$.\par
$(A)$ For $\mathbb{U}_0\in\mathbb{W}_2^{1,0}\cap\mathbb{L}_q$ and $t\geqq1$, we have the following estimates:
\begin{align*}
&\|\mathbb{U}(t)\|_2\leqq Ct^{-\sigma}(\|\mathbb{U}_0\|_q+\|\mathbb{U}_0\|_{\mathbb{W}_2^{1,0}}),\quad\sigma=\frac{3}{2}\left(\frac{1}{q}-\frac{1}{2}\right),\\
&\|\partial_x\mathbb{U}(t)\|_2+\|\partial_t\mathbb{U}(t)\|_2\leqq Ct^{-\sigma-\frac{1}{2}}(\|\mathbb{U}_0\|_q+\|\mathbb{U}_0\|_{\mathbb{W}_2^{1,0}}),\\
&\|\partial_x^2\mathbb{PU}(t)\|_2\leqq Ct^{-\frac{3}{2q}}(\|\mathbb{U}_0\|_q+\|\mathbb{U}_0\|_{\mathbb{W}_2^{1,0}}),\\
&\|\mathbb{U}(t)\|_{\mathbb{W}_\infty^{0,1}}\leqq Ct^{-\frac{3}{2q}}(\|\mathbb{U}_0\|_q+\|\mathbb{U}_0\|_{2,2}).
\end{align*}\par
$(B)$ If $\mathbb{U}_0\in\mathbb{W}_2^2\cap\mathbb{L}_q$ and $\mathbb{P}\mathbb{U}_0=0$ on $\partial \Omega $, then we have the following estimates for $t\geqq1$:
\begin{align*}
&\|\partial_x^2(\mathbb{I}-\mathbb{P})\mathbb{U}(t)\|_2+\|\partial_x^3\mathbb{P}\mathbb{U}(t)\|_2+\|\partial_x\partial_t\mathbb{U}(t)\|_2 \leqq Ct^{-\frac{3}{2q}}(\|\mathbb{U}_0\|_q+\|\mathbb{U}_0\|_{2,2}),\\
&\|\mathbb{U}(t)\|_{\mathbb{W}_\infty^{1,2}}\leqq Ct^{-\frac{3}{2q}}(\|\mathbb{U}_0\|_q+\|\mathbb{U}_0\|_{\mathbb{W}_2^{2,1}}).
\end{align*}\par
$(C)$ If we assume that $q>1$ additionally, then for $\mathbb{U}_0\in\mathbb{W}_2^{2,1}\cap\mathbb{W}_q^{1,0}$ and $t\geqq1$ we have the following estimates:
\begin{align*}
&\|\partial_x^2(\mathbb{I}-\mathbb{P})\mathbb{U}(t)\|_2+\|\partial_x^3\mathbb{P}\mathbb{U}(t)\|_2+\|\partial_x\partial_t\mathbb{U}(t)\|_2\leqq Ct^{-\frac{3}{2q}}(\|\mathbb{U}_0\|_{\mathbb{W}_q^{1,0}}+\|\mathbb{U}_0\|_{\mathbb{W}_2^{2,1}}),\\
&\|\mathbb{U}(t)\|_{\mathbb{W}_\infty^{1,2}}\leqq Ct^{-\frac{3}{2q}}(\|\mathbb{U}_0\|_{\mathbb{W}_q^{1,0}}+\|\mathbb{U}_0\|_{\mathbb{W}_2^{2,1}}).
\end{align*}
\end{Theorem}
Theorem \ref{the1.3} will be proved by combination of the $L_2$ type estimate in the $\mathbb{R}^3$ case and the local energy decay estimate of $\{e^{t\mathbb{A}}\}_{t\geq0}$, via a cut-off technique in Sect.4. Such a cut-off technique is already known as an effective method in the study of the $L_2$ type estimate for the semigroups in the exterior domain.
\section{Some Properties of the Semigroup $\{e^{t\mathbb{A}}\}_{t\geq0}$}
In this section, we shall summarize some properties of the semigroup $\{e^{t\mathbb{A}}\}_{t\geq0}$, which was obtained by Kobayashi \cite{Kobayashi1997}. Let $\rho(\mathbb{A})$ be the resolvent set of the operator $\mathbb{A}$. Then, we have the following lemma.
\begin{Lemma}\label{lem2.1}
$\mathbb{A}$ is a closed linear operator in $\mathbb{W}_2^{1,0}(\Omega)$ and
$$\rho(\mathbb{A})\supset\Sigma=\{\kappa\in\mathbb{C}\mid C\ \mathrm{Re}\ \kappa+(\mathrm{Im}\ \kappa)^2>0\},$$
where $C$ is a constant depending only on $\mu,\lambda$. Moreover, the following properties are valid: There exist positive constants $\kappa_0$ and $\delta<\frac{\pi}{2}$ such that for any $\kappa_0\in\Sigma_\delta=\{\kappa\in\mathbb{C}\mid|arg\kappa\mid\leqq\pi-\delta\},$
\begin{align*}
\begin{aligned}
&|\kappa|\|(\kappa-\mathbb{A})^{-1}\mathbb{F}\|_{\mathbb{W}_2^{1,0}}+\|\mathbb{P}(\kappa-\mathbb{A})^{-1}\mathbb{F}\|_{2,2}\leqq C(\kappa_0,\delta)\|\mathbb{F}\|_{\mathbb{W}_2^{1,0}},\quad
\text{if}\ \mathbb{F}\in\mathbb{W}_2^{1,0};\\
&|\kappa|^{\frac{1}{2}}\|(\mathbb{I}-\mathbb{P})(\kappa-\mathbb{A})^{-1}\mathbb{F}\|_{2,2}+|\kappa|^{-\frac{1}{2}}\|\mathbb{P}(\kappa-\mathbb{A})^{-1}\mathbb{F}\|_{3,2}\leqq C(\kappa_0,\delta)\|\mathbb{F}\|_{\mathbb{W}_2^{2,1}},\quad
\text{if}\ \mathbb{F}\in\mathbb{W}_2^{2,1}.\end{aligned}
\end{align*}
\end{Lemma}
By Lemma\ref{lem2.1} we see that $\mathbb{A}$ generates an analytic semigroup $\{e^{t\mathbb{A}}\}_{t\geqq0}$ on $\mathbb{W}_2^{1,0}$ and moreover we have the following estimates for $0<t\leqq2$:
\begin{align}\label{3.1.}
\|(-\mathbb{A})^ke^{t\mathbb{A}}\mathbb{U}\|_{\mathbb{W}_2^{1,0}}&\leqq C(p,k)t^{-k}\|\mathbb{U}\|_{\mathbb{W}_2^{1,0}},\quad
\text{for}\ \mathbb{U}\in\mathbb{W}_2^{1,0},\ k\geqq0\ \mathrm{integer},
\end{align}
\begin{align}\label{2.2}
\|e^{t\mathbb{A}}\mathbb{U}\|_{1,2,\Omega}&\leqq C(p)t^{-\frac{1}{2}}\|\mathbb{U}\|_{\mathbb{W}_2^{1,0}},\quad\mathrm{for~}\mathbb{U}\in\mathbb{W}_2^{1,0},
\end{align}
\begin{align}\label{2.3}
\|(\mathbb{I}-\mathbb{P})e^{t\mathbb{A}}\mathbb{U}\|_{2,2,\Omega}&\leqq C(p)t^{-\frac{1}{2}}\|\mathbb{U}\|_{\mathbb{W}_2^{2,1}},\quad\mathrm{for~}\mathbb{U}\in\mathbb{W}_2^{2,1},
\end{align}
\begin{align}\label{2.4}
\|\mathbb{P}e^{t\mathbb{A}}\mathbb{U}\|_{3,2,\Omega}&\leqq C(p)t^{-\frac{3}{2}}\|\mathbb{U}\|_{\mathbb{W}_2^{2,1}},\quad\mathrm{for~}\mathbb{U}\in\mathbb{W}_2^{2,1}.
\end{align}
By Lemma\ref{lem 7.3} in the Appendix below, we have
\begin{align}\label{3.5.}
C_1\|\mathbb{U}\|_{\mathbb{W}_2^{1,2}}\leqq\|\mathbb{U}\|_{\mathbb{W}_2^{1,0}}+\|-\mathbb{A}\mathbb{U}\|_{\mathbb{W}_2^{1,0}}\leqq C_2\|\mathbb{U}\|_{\mathbb{W}_2^{1,2}}
\end{align}
for suitable constants $C_1$ and $C_2$ provides that $\mathbb{U}\in\mathcal{D}_2(\mathbb{A})$, $\mathcal{D}_2(\mathbb{A})=\{\mathbb{U}=(n,\mathbf{v},E)\in\mathbb{W}_2^{1,2}\mid\mathbb{P}\mathbb{U}|_{\partial\Omega}=0\}$. If $\mathbb{U}\in\mathcal{D}_2(\mathbb{A})$, then $-\mathbb{A}e^{t\mathbb{A}}\mathbb{U}=e^{t\mathbb{A}}(-\mathbb{A}\mathbb{U})$, and therefore by \eqref{3.1.} and \eqref{3.5.}
\begin{align}\label{2.6}
\|e^{t\mathbb{A}}\mathbb{U}\|_{\mathbb{W}_2^{1,2}}\leqq C\|\mathbb{U}\|_{\mathbb{W}_2^{1,2}}.
\end{align}
The next lemma is concerned with a local decay property of $\{e^{t\mathbb{A}}\}_{t\geq0}$, which was obtained by Ishigaki and Kobayashi\cite{Ishigaki-Kobayashi2025}. Set
\begin{align*}
\mathbb{W}_{2,b}^{k,m}&=\{\mathbb{U}\in\mathbb{W}_2^{k,m}\mid\mathbb{U}(x)=0\mathrm{~for}|x|\geqq b\},\\
\Omega_{b}&=\Omega\cap B_b,\quad B_b=\{x\in\mathbb{R}^3\mid|x|<b\}.
\end{align*}
\begin{Lemma}\label{lem2.2}
(Local energy decay). Let $b_0$ be a fixed number such that $B_{b_0}\supset\mathbb{R}^3\setminus\Omega.$ Suppose that $b>b_0$ and $\Omega_b=\Omega\cap B_b.$ If $\mathbb{U}_0=(n_0,v_0,E_0)$ satisfies $\mathbb{U}_0\in \mathbb{W}_2^{1,0}(\Omega),\nabla n _0+\mathrm{div} {E_0}^T=0,\nabla E_0\in W^{2}_2(\Omega)\cap \mathcal{W}_2^1(D)(\Omega)$ and $\mathbb{U}_0=0 (x\notin\Omega_b)$, then the following local energy decay estimates of $\mathbb{U}(t)=\mathbb{T}(t)\mathbb{U}_0$ hold for all $\begin{aligned}m=0,1,2,\ldots\end{aligned}$ and $t\geqq1$:
\begin{equation*}
\|\partial_t^m\mathbb{T}(t)\mathbb{U}_0\|_{\mathbb{W}_2^{1,2}(\Omega_b)}\leqq C_{b,b_0,m}t^{-2-m}\|\mathbb{U}_0\|_{\mathbb{W}_2^{1,0}(\Omega)},
\end{equation*}
where $C_{b,b_0,m}$ is a constant depending only on $b,b_0,m$.
\end{Lemma}
\section{$L_2$ Estimates of Solutions to the Linearized Equations in $\mathbb{R}^3$}
In this section, we shall show $L_2$ estimate of solutions to the initial value problem in $\mathbb{R}^3$. For this purpose, we need to consider the following linearized system
\begin{align}\label{3.1'}
\mathbb{U}_t=\mathbb{AU}\ \text{in}\ [0,\infty)\times\mathbb{R}^3,\quad
\mathbb{U}(0)=\mathbb{F}\ \text{in}\ \mathbb{R}^3,
\end{align}
where $\mathbb{U}=(n,v,E)$ and $\mathbb{F}=(f_1,f_2,f_3)$. Define
\begin{align}\label{3.2'}
\mathbb{U}(t)=\mathbb{E}(t)\mathbb{F}.
\end{align}
\begin{Theorem}\label{th3.1}
Let $\mathbb{E}(t)$ be the solution operator of \eqref{3.1'} defined by \eqref{3.2'}. Then, we have the decomposition:
$$\mathbb{E}(t)\mathbb{F}=\mathbb{E}_0(t)\mathbb{F}+\mathbb{E}_\infty(t)\mathbb{F},$$
where $\mathbb{E}_0(t)$ and $\mathbb{E}_\infty(t)$ have the following properties:\par
$(1)\forall\ell,m\geqq0$ integers,
	$$|\partial_t^m\partial_x^\ell\mathbb{E}_0(t)\mathbb{F}|_2\leqq C(m,\ell,q)t^{-\frac{3}{2}(\frac{1}{q}-\frac{1}{2})-\frac{m+\ell}{2}}|\mathbb{F}|_q,\quad \forall t\geqq1,$$
where  $1\leqq q\leqq2$, and
$$|\partial_t^m\partial_x^\ell\mathbb{E}_0(t)\mathbb{F}|_2\leqq C(m,\ell,p)|\mathbb{F}|_q,\quad 0<\forall t\leqq2,$$
where $1\leq q\leq2$.\par
$(2) \text{Set} (\ell)^+=\ell\mathrm{~if~}\ell\geq0\mathrm{~and~}(\ell)^+=0\mathrm{~if~}\ell<0.$ There exists a $c>0$ such that for $\forall l,m,n \geqq 0$ integers,
$$\begin{gathered}\begin{aligned}|\partial_t^m\partial_x^\ell(\mathbb{I}-\mathbb{P})\mathbb{E}_\infty(t)\mathbb{F}|_2\end{aligned} \leqq C(m,\ell,n)e^{-ct}\left[t^{-\frac{n}{2}}|\mathbb{F}|_{(2m+\ell-n-1)^{+},2}+|\mathbb{F}|_{\mathbb{W}_{2}^{\ell,(\ell-1)^{+}}}\right] ;\\
\begin{aligned}|\partial_t^m\partial_x^\ell\mathbb{P}\mathbb{E}_\infty(t)\mathbb{F}|_2\end{aligned} \leqq C(m,\ell,n)e^{-ct}\left[t^{-\frac{n}{2}}|\mathbb{F}|_{W_2^{(2m+\ell-n-1)^+,(2m+\ell-n)^+}}+|\mathbb{F}|_{W_2^{(\ell-1)^+,(\ell-2)^+}}\right],  \quad  \forall t>0\end{gathered}$$
$(3)$  Let $c$ be the same as in $(2)$. Then, we have
\begin{align*}
&|\partial_t^m\partial_x^\ell(\mathbb{I}-\mathbb{P})\mathbb{E}_\infty(t)\mathbb{F}|_\infty \leqq C(m,\ell,n)t^{-\left(\frac{n}{2}+\frac{3}{4}\right)}e^{-ct}|\mathbb{F}|_{(2m+\ell-n-1)^+,2}+C(m,\ell)e^{-ct}\left[|\mathbb{F}|_{\mathbb{W}_2^{(\ell,\ell+1)}}+|(\mathbb{I}-\mathbb{P})\mathbb{F}|_{H^{l+2}}\right]\\
&|\partial_t^m\partial_x^\ell\mathbb{P}\mathbb{E}_\infty(t)\mathbb{F}|_\infty \leqq C(m,\ell,n)e^{-ct}\left[t^{-\left(\frac{n}{2}+\frac{3}{4}\right)}|\mathbb{F}|_{\mathbb{W}_2^{(2m+\ell-n-1)^+,(2m+\ell-n)^+}} +\left|\mathbb{F}\right|_{\mathbb{W}_2^{\ell+1,\ell}}\right],\quad\forall t>0.
\end{align*}
\end{Theorem}
To simplify the calculation of eigenvalues in spectral analysis, we make some transformation. The strategy follows Hu and Wu \cite{Hu-Wu2013}. For
\begin{align*}
\frac{\partial^2(E^{ij})}{\partial x_i\partial x_j}
&=\mathrm{div}\mathrm{div}(E^T)\\
&=\mathrm{div}\mathrm{div}[(1+n)(E+I)^T]-\mathrm{div}\mathrm{div}(nI+E^T)\\
&=-\triangle n-\mathrm{div}\mathrm{div}(nE),
\end{align*}
by applying $\mathrm{div}$ to $\eqref{2.8.}_2$, we have
\begin{align}\label{2.25}
(\mathrm{div}v)_t-(2\mu+\lambda)\triangle \mathrm{div}v+2\triangle n=0.
\end{align}
Next, taking the transpose of $\eqref{2.8.}_3$ and then minusing $\eqref{2.8.}_3$, we have
\begin{align}\label{2.28}
(E^T-E)_t+\mathcal{L}=h^T-h-v\cdot\nabla(E^T-E).
\end{align}
where $\mathcal{L}=\nabla v-(\nabla v)^T=\mathrm{curl}v.$ Note the condition $F^{lk}\nabla_lF^{ij}=F^{lj}\nabla_lF^{ik}$ for all $t\geqq0$, which means that
\begin{align}
\nabla_kE^{ij}+E^{lk}\nabla_lE^{ij}=\nabla_jE^{ik}+E^{lj}\nabla_lE^{ik},\quad\forall t\geq0.
\end{align}
Thus we have
\begin{align*}
\nabla_j\nabla_kE^{ik}-\nabla_i\nabla_kE^{jk}
&=\nabla_k\nabla_jE^{ik}-\nabla_k\nabla_iE^{jk}\\
&=\nabla_k\nabla_kE^{ij}-\nabla_k\nabla_kE^{ji}+\nabla_k(E^{lk}\nabla_lE^{ij}-E^{lj}\nabla_lE^{ik})\\
&-\nabla_k(E^{lk}\nabla_lE^{ji}-E^{li}\nabla_lE^{jk}).
\end{align*}
Thus by applying $\mathrm{curl}$ to $\eqref{2.8.}_2$, we have
\begin{align}\label{2.31}
\mathcal{L}_t-\mu\triangle\mathcal{L}+\triangle(E^T-E)=0.
\end{align}
Notice that the system \eqref{2.28}-\eqref{2.31} takes an argument similar to the system $\eqref{2.8.}_2$ and \eqref{2.25}. We find that the linearized system$\eqref{2.8.}_1$  and \eqref{2.25} depends on $(n,\mathrm{div}v)$ while the linearized system \eqref{2.28}-\eqref{2.31} depends only on $(\mathcal{L},E^T-E)$. Denote by $\Lambda^{s}$ the pseudodifferential operator defined by
$$\Lambda^su=\mathcal{F}^{-1}(|\xi|^s\hat{u}(\xi)),$$
and let
$$d=\Lambda^{-1}\mathrm{div}v$$
be the compressible part of the velocity, and
$$\omega=\Lambda^{-1}\mathcal{W}=\Lambda^{-1}\mathrm{curl}v$$
be the incompressible part of the velocity. We thus obtain
\begin{align}\label{3.1}
\begin{cases}
n_t+\Lambda d=0,\\
d_t-(2\mu+\lambda)\triangle d-2\Lambda n=0,
\end{cases}
\end{align}
and
\begin{align}\label{3.2}
\begin{cases}
(E^T-E)_t+\Lambda\omega=0,\\
\omega_t-\mu\triangle\omega-\Lambda(E^T-E)=0
\end{cases}
\end{align}
Indeed, as the definition of $d$ and $\omega$, and the relation
$$v=-\Lambda^{-1}\nabla d-\Lambda^{-1}\mathrm{div}\omega$$
involve pseudodifferential operators of degree zero, the estimate in the space $H^l(\mathbb{R}^3)$ for the original function $v$ will be the same as for $(d,\omega)$. Then,
Now, we discuss the system \eqref{3.1}. Defined  $B$ be a matrix-valued differential operator given by
$$B:=\begin{pmatrix}0&-\Lambda\\2\Lambda&(2\mu+\lambda)\triangle\end{pmatrix}.$$\par
Now we aim to analyze the differential operator $B$ in terms of its Fourier expression $\hat{B}(\xi)$ and to show the long time properties of the semigroup $e^{tB}$ For this purpose, we need to consider the following linearized system
\begin{align}\label{3.1.1}
U_t-BU=0\ \text{in}\ [0,\infty)\ \times\mathbb{R}^3,\quad U(0)=G\ \text{in}\ \mathbb{R}^3,
\end{align}
where $U=(n,d)$ and $G=(g_1,g_2)$. Applying the Fourier transform to system \eqref{3.1.1}, we have
\begin{align}\label{3.2.2}
U(t)=K(t)G=\mathcal{F}^{-1}(e^{t\hat{B}(\xi)}\hat{G}(\xi)),
\end{align}
 and $\hat{B}(\xi)$ is defined as
$$\hat{B}(\xi):=\begin{pmatrix}0&-|\xi|\\2|\xi|&-(2\mu+\lambda)|\xi|^2\end{pmatrix}.$$
The characteristic polynomial of $\hat{B}(\xi)$ is $\begin{aligned}\kappa^2+(2\mu+\lambda)|\xi|^2\kappa+2|\xi|^2\end{aligned}$, which implies the eigenvalues are
$$\kappa_\pm=-\left(\mu+\frac{1}{2}\lambda\right)|\xi|^2\pm\frac{1}{2}i\sqrt{8|\xi|^2-(2\mu+\lambda)^2|\xi|^4}.$$
The semigroup $e^{t\hat{B}(\xi)}$ is expressed as
\begin{align*}e^{t\hat{B}(\xi)}&=e^{\kappa_+t}\frac{\hat{B}(\xi)-\kappa_-I}{\kappa_+-\kappa_-}+e^{\kappa_-t}\frac{\hat{B}(\xi)-\kappa_+I}{\kappa_--\kappa_+}\\
&=\begin{pmatrix}\frac{\kappa_+e^{\kappa_-t}-\kappa_-e^{\kappa_+t}}{\kappa_+-\kappa_-}&-\frac{|\xi|(e^{\kappa_+t}-e^{\kappa_-t})}{\kappa_+-\kappa_-}\\\frac{2|\xi|(e^{\kappa_+t}-e^{\kappa_-t})}{\kappa_+-\kappa_-}&\frac{\kappa_+e^{\kappa_-t}-\kappa_-e^{\kappa_+t}}{\kappa_+-\kappa_-}-\frac{(2\mu+\lambda)|\xi|^2(e^{\kappa_+t}-e^{\kappa_-t})}{\kappa_+-\kappa_-}\end{pmatrix}.
\end{align*}
When $|\xi|\leqq\frac{2\sqrt{2}}{2\mu+\lambda}$,
$$\kappa_{+}=\overline{\kappa_{-}}=-(\mu+\frac{1}{2}\lambda)|\xi|^{2}+\frac{1}{2}i\sqrt{8|\xi|^{2}-(2\mu+\lambda)^{2}|\xi|^{4}}.$$
and there exists a positive constant $ r_{1}\leqq\frac{2}{2\mu+\lambda}$, such that $\kappa_{+},\kappa_{-}$ has a Taylor series expansion for $|\xi|<r_1$ as follows:
$$\kappa_{+}=\overline{\kappa_{-}}=-(\mu+\frac{1}{2}\lambda)|\xi|^{2}+\sqrt{2}|\xi|i+O(|\xi|^{3}).$$
Similarly, when $|\xi|\geqq\frac{2\sqrt{2}}{2\mu+\lambda}$,
$$\kappa_{\pm}=-(\mu+\frac{1}{2}\lambda)|\xi|^{2}\pm\frac{1}{2}\sqrt{-8|\xi|^{2}+(2\mu+\lambda)^{2}|\xi|^{4}}.$$
and there exists a positive constant $ r_{2}\geqq\frac{2}{2\mu+\lambda}$, such that $\kappa_{+},\kappa_{-}$ has a Laurent series expansion for $|\xi|>r_2$ as follows:
$$\kappa_+(\xi)=-\frac{2}{2\mu+\lambda}-\frac{4}{(2\mu+\lambda)^{3}}\frac{1}{|\xi|^{2}}+O(\frac{1}{|\xi|^{4}}),$$
$$\kappa_-(\xi)=-(2\mu+\lambda)|\xi|^{2}+\frac{2}{2\mu+\lambda}+\frac{4}{(2\mu+\lambda)^{3}}\frac{1}{|\xi|^{2}}+O(\frac{1}{|\xi|^{4}}).$$
The matrix exponential has the spectral resolution
\begin{align*}
e^{t\hat{B}(\xi)}&=e^{t\kappa_+(\xi)}\mathbb{P}_+(\xi)+e^{t\kappa_-(\xi)}\mathbb{P}_-(\xi)\\
&=e^{t\kappa_+}\frac{\hat{B}(\xi)-\kappa_-I}{\kappa_+-\kappa_-}+e^{t\kappa_-}\frac{\hat{B}(\xi)-\kappa_+I}{\kappa_--\kappa_+}.
\end{align*}
Next, we prove Theorem \ref{th3.1}. The strategy follows Kobayashi and Shibata \cite{Kobayashi-Shibata1999}.
\begin{proof}
Let $\varphi_{0}(\xi)$ be a function in $C_0^\infty(\mathbb{R}^3)$ such that $\varphi_{0}(\xi)=1$ for $|\xi|\leqq r_1/2$ and $\varphi_{0}(\xi)=0$  for $|\xi|\geqq r_1$ and set
$$\Psi_0(t)G(x)=\sum_{j={\pm}}\mathcal{F}^{-1}[e^{t\kappa_j(\xi)}\varphi_0(\xi)\mathbb{P}_j(\xi)\hat{G}(\xi)](x).$$
in order to estimate $\Psi_0(t)G(x)$, we consider the function:
$$H_\ell(t)f(x)=\mathcal{F}^{-1}[e^{\kappa(\xi)t}h_\ell(\xi)\varphi_0(\xi)\hat{g}(\xi)](x),$$
where $\kappa(\xi)$ and $h_\ell(\xi)$ satisfy the condition:
$$-\beta_{0}|\xi|^{2}\leqq\mathrm{Re}\kappa(\xi)\leqq-\beta_{1}|\xi|^{2},\quad|\xi^{\alpha}\partial_{\xi}^{\alpha}\kappa(\xi)|\leqq C(\alpha)|\xi|\mathrm{~and~}|\xi^{\alpha}\partial_{\xi}^{\alpha}h_{\ell}(\xi)|\leqq C(\alpha)|\xi|^{\ell}$$
for any $\xi\in\mathbb{R}^3\backslash\{0\}$ with $|\xi|\leqq r_1$.
First, we consider the case of $t\geqq0$ and $1\leqq q\leqq2\leqq p\leqq\infty$. Then,
$$F_\ell(t)g(x)=\left[\mathcal{F}^{-1}[e^{-\beta_1|\xi|^2t/4}]*\mathcal{F}^{-1}[e^{\mu(\xi)t}h_\ell(\xi)\varphi_0(\xi)]*\mathcal{F}^{-1}[e^{-\beta_1|\xi|^2t/4}]*g\right](x)$$
and
$$\mu(\xi)=\kappa(\xi)+\beta_1|\xi|^2/2,$$
where $*$ means the convolution and $\mathrm{Re}\mu(\xi)\leqq-\beta_1|\xi|^2/2.$ Since
$$\mathcal{F}^{-1}[e^{-\beta_1|\xi|^2t/4}]=(\beta_1\pi t)^{-3/2}e^{-|x|^2/\beta_1t},$$
we obtain after using Young's inequality that
\begin{align}\label{3.6}
\|\mathcal{F}^{-1}[e^{-\beta_1|\xi|^2{t/4}}]*u\|_p\leqq Ct^{-\frac{3}{2}\left(\frac{1}{q}-\frac{1}{p}\right)}\|u\|_q,
\end{align}
where $1\leqq q\leqq p\leqq\infty$. By Parseval's formula, one has
$$\begin{aligned}\|\mathcal{F}^{-1}[e^{\mu(\xi)t}h_\ell(\xi)\varphi_0(\xi)]*u\|_2\leqq C\left\{\int_{\mathbb{R}^3}e^{-\beta_1|\xi|^2t}|\xi|^{2\ell}|\hat{u}(\xi)|^2d\xi\right\}^{\frac{1}{2}}.
\end{aligned}$$
Let $\xi=\frac{\eta}{\sqrt{t}}$, we thus get
$$\begin{aligned}\|\mathcal{F}^{-1}[e^{\mu(\xi)t}h_\ell(\xi)\varphi_0(\xi)]*u\|_2&\leqq
C\int_{\mathbb{R}^{3}}e^{-\beta_1\eta^{2}}
\frac{\left|\eta\right|^{2\ell}}{t^{\ell}}|\hat{h}(\frac{\eta}{\sqrt{t}})|^{2}t^{-\frac{3}{2}}d\eta\\&
\leq Ct^{-\frac{\ell}{2}}\|u\|_2.\end{aligned}$$
Therefore, for $t \geqq1$, $s \leqq2$ and $1\leqq s\leqq2\leqq p\leqq\infty$, we have
\begin{align*}
\|H_\ell(t)g\|_p&=\left\|[\mathcal{F}^{-1}[e^{-\beta_{1}|\xi|^{2}t/4}]*\mathcal{F}^{-1}[e^{\mu(\xi)t}g_{\ell}(\xi)\varphi_{0}(\xi)]*\mathcal{F}^{-1}[e^{-\beta_{1}|\xi|^{2}t/4}]*g](x)\right\|_p \\ &
&=Ct^{-\frac{3}{2}(\frac{1}{q}-\frac{1}{p})-\frac{\ell}{2}}{\left\|g\right\|_q},
\end{align*}
or
\begin{align*}
\|H_\ell(t)g\|_p&=\left\|[\mathcal{F}^{-1}[e^{-\beta_{1}|\xi|^{2}t/4}]*\mathcal{F}^{-1}[e^{\mu(\xi)t}g_{\ell}(\xi)\varphi_{0}(\xi)]*\mathcal{F}^{-1}[e^{-\beta_{1}|\xi|^{2}t/4}]*g](x)\right\|_p \\
&\leqq Ct^{-\frac{3}{2}(\frac{1}{2}-\frac{1}{p})-\frac{\ell}{2}-\frac{3}{2}(\frac{1}{s}-\frac{1}{2})}\left\|\hat{g}\right\|_\frac{s}{s-1}.
\end{align*}
Next, we consider the case of $0<t\leqq1$, by the Fourier multiplier theorem we have
$$\left\|\partial_x^\alpha H_\ell(t)g\right\|_q\leqq C(\alpha,p)\left\|g\right\|_q\quad0<\forall t\leqq1,\quad\forall\alpha,$$
for $1<q<\infty$, together with Sobolev's imbedding theorem implies that
$$\left\|H_\ell(t)g\right\|_p\leqq C(p,q)\left\|g\right\|_q\quad0<\forall t\leqq1,$$
for $1<q\leqq p\leqq \infty$ and $1<q<\infty$. If we choose $\psi(\xi)\in C_0^\infty(\mathbb{R}^n)$ so that $\varphi_0(\xi)=1$ on the supp $\psi(\xi)$, then we have
$$H_\ell(t)f=(H_\ell(t)\mathcal{F}^{-1}[\psi])*g,$$
which combined with Young's inequality gives
$$\left\|H_\ell(t)g\|_p\leqq\|H_\ell(t)\mathcal{F}^{-1}[\psi]\right\|_p\left\|g\right\|_1\leqq C\left\|g\right\|_1$$
for $1<p<\infty$. Finally, we have
$$\|H_\ell(t)g\|_\infty\leqq\int_{\mathbb{R}^3}\left|e^{\kappa(\xi)t}h_\ell(\xi)\varphi_0(\xi)\hat{g}(\xi)\right|d\xi\leqq C\|g\|_1.$$
Combining these inequalities, we have
$$\|G_\ell(t)f\|_p\leqq C\|f\|_q,\quad0<t\leqq1,$$
for $1\leqq q\leqq p \leqq \infty$ and $(p,q)\neq(1,1),(\infty,\infty).$\par
Applying these estimates with $h_\ell(\xi)\hat{g}(\xi)=\xi^\alpha\kappa(\xi)^m\mathbb{P}_j(\xi)\hat{G}(\xi)$ with $|\alpha|=\ell$ and $\kappa(\xi)=\kappa_j(\xi)$, immediately we have
$$\begin{aligned}&\|\partial_t^m\partial_x^\ell\Psi_0(t)G\|_p\leqq Ct^{-\frac{3}{2}\left(\frac{1}{q}-\frac{1}{p}\right)-\frac{(m+\ell)}{2}}\|G\|_q,&\mathrm{for~}t\geqq1,\quad1\leqq q\leqq2\leqq p\leqq\infty,\\&\|\partial_t^m\partial_x^\ell\Psi_0(t)G\|_p\leqq\|G\|_q,&\mathrm{for~}\begin{cases}0<t\leqq1,\quad1\leqq q\leqq p\leqq\infty,\\(p,q)\neq(1,1),(\infty,\infty).&\end{cases}\end{aligned}$$\par
Next, we estimate the term corresponding to the large $|\xi|$ and just consider the $L_2$ norm. Set
$$\Psi_\infty(t)G=\sum_{j={\pm}}\mathcal{F}^{-1}\left[e^{\kappa_j(\xi)t}\varphi_\infty(\xi)\mathbb{P}_j(\xi)\hat{G}(\xi)\right](x),$$
where $\varphi_{\infty}(\xi)$ is a function in $C^\infty(\mathbb{R}^3)$ such that $\varphi_{\infty}(\xi)=1$ for $|\xi|\geqq 3 r_2$ and $\varphi_{\infty}(\xi)=0$ for $|\xi|\leqq 2 r_2$. In order to estimate $\Psi_\infty(t)G$, we shall consider the function:
$$[R_\ell^j(t)g](x)=\mathcal{F}^{-1}\left[e^{\kappa_j(\xi)t}h_\ell(\xi)\varphi_\infty(\xi)\hat{g}(\xi)\right](x),\quad j={\pm},$$
where $g_\ell(\xi)$ satisfy the condition:
$$|\xi^\alpha\partial_\xi^\alpha h_\ell(\xi)|\leqq C(\alpha)|\xi|^{-\ell},\quad\forall\xi\in\mathbb{R}^3.$$
For $\kappa_{-}$, there exist $\beta_0>\beta_1>0$ such that
$$-\beta_{0}|\xi|^{2}\leqq\mathrm{Re}\kappa_{-}(\xi)\leqq-\beta_{1}|\xi|^{2}
\quad \mathrm{and}\quad
|\xi^{\alpha}\partial_{\xi}^{\alpha}\kappa_{-}(\xi)|\leqq C(\alpha)|\xi|^{2}$$
for $|\xi|\geqq r_2$. Therefore, by induction we have
$$|\xi^\alpha\partial_\xi^\alpha e^{\kappa_-(\xi)t}|\leqq C(\alpha)e^{-\frac{\beta_1}{2^{|\alpha|}}|\xi|^2t}$$
for $|\xi|\geqq r_2$ and $t>0$, which implies that
$$\left|\xi^{\alpha}\partial_{\xi}^{\alpha}\left[\frac{e^{\kappa_-(\xi)t}\xi^{\beta}\kappa_-(\xi)^{m}h_{\ell}(\xi)\varphi_{\infty}(\xi)}{(1+|\xi|^{2})^{(\frac{2m+|\beta|-\ell-n}{2})^{+}}}\right]\right|\leqq C(\alpha,\beta,m,\ell,n)t^{-\frac{n}{2}}e^{-\frac{\beta_{1}}{4^{|\alpha|}}|\xi|^{2}t}$$
for $|\xi|\geqq r_2$ and $t>0$. Since
$$\|\mathcal{F}^{-1}[(1+|\xi|^2)^{\frac{a}{2}}\hat{g}(\xi)]\|_2\leqq C(a,p)\|g\|_{a,2},$$
we get from the Fourier multiplier theorem that
\begin{align}\label{3.7}
\|\partial_t^m\partial_x^\beta R_\ell^{+}(t)f\|_2\leqq C(m,\beta,n,p)t^{-\frac{n}{2}}e^{-ct}\|g\|_{(2m+|\beta|-\ell-n)^+,2}\quad\forall t>0,
\end{align}
where $c=\frac{\beta_1r_2^2}{16}$.\par
Next, we shall consider the case that $j=+$. There exist constants $\beta_0>\beta_1>0$ such that
$$-\beta_0\leqq\mathrm{Re}\kappa_+(\xi)\leqq-\beta_1$$
for $|\xi|\geqq r_2$.
There exist a positive constant $\gamma$ and a function $\mu_+(\xi)$ such that
$$\kappa_+(\xi)=-\gamma+\mu_+(\xi)\mathrm{~and~}|\xi^\alpha\partial_\xi^\alpha\mu_+(\xi)|\leqq C(\alpha)|\xi|^{-2}$$
for $|\xi|\geqq r_2$. Since
$$\left|\xi^\alpha\partial_\xi^\alpha\left[\frac{e^{\kappa_+(\xi)t}\xi^\beta\kappa_+(\xi)^mh_\ell(\xi)\varphi_\infty(\xi)}{(1+|\xi|^2)^{(\frac{|\beta|-\ell}{2})^+}}\right]\right|\leqq C(\alpha,\beta,m,\ell)(1+t)^{|\alpha|}e^{-\beta_1t}\quad \forall t>0,$$
by the Fourier multiplier theorem, we obtain
\begin{align}\label{3.9}
\|\partial_t^m\partial_x^\beta R_\ell^+(t)f\|_2&\leqq C(m,\beta,\ell)(1+t)^3e^{-\beta_1t}\|g\|_{(|\beta|-\ell)^+,2} \notag \\&\leqq C(m,\beta,\ell,\beta_1)e^{-\frac{\beta_1}{2}t}\|g\|_{(|\beta|-\ell)^+,2}
\end{align}
for any $t>0$.\par
Now, using \eqref{3.7} and \eqref{3.9}, we shall estimate $\Psi_\infty(t)G$. Put
$$\sum_{j={\pm}}e^{\kappa_j(\xi)t}\varphi_\infty(\xi)\mathbb{P}_j(\xi)\hat{G}(\xi)=(\hat{n}(t,\xi),\hat{d}(t,\xi)).$$
We write
$$\hat{n}(t,\xi)=\sum_{k={\pm}}\sum_{j=1,2}e^{\kappa_k(\xi)t}A_{kj}^1(\xi)\varphi_\infty(\xi)\hat{g}_j(\xi),$$
where
\begin{align}\label{3.10}
&A_{+1}^1(\xi)=1+p_{+1}^1(\xi)\mathrm{~and~}|\xi^\alpha\partial_\xi^\alpha p_{+1}^1(\xi)|\leqq C(\alpha)|\xi|^{-2},\\ \notag
&\left|\xi^\alpha\partial_\xi^\alpha A_{k2}^1(\xi)\right|\leqq C_\alpha|\xi|^{-1}\mathrm{~and~}\left|\xi^\alpha\partial_\xi^\alpha A_{-2}^1(\xi)\right|\leqq C_\alpha|\xi|^{-2}
\end{align}
for $|\xi|\geqq r_2$ and $\kappa={\pm}$. By Fourier mulitiplier theorem and Sobolev's inequality, we have
$$\begin{aligned}\left\|\partial_t^m\partial_x^\ell\mathcal{F}^{-1}\left[e^{\kappa_+(\xi)t}A_{+1}^1(\xi)\varphi_\infty(\xi)\hat{g}_1(\xi)\right]\right\|_2&\leqq C(m,\ell)e^{-\gamma t}(1+t)^m\|g_1\|_{\ell,2},\\ \left\|\partial_t^m\partial_x^\ell\mathcal{F}^{-1}\left[e^{\kappa_+(\xi)t}A_{+1}^1(\xi)\varphi_\infty(\xi)\hat{g}_1(\xi)\right]\right\|_\infty&\leqq C(m,\ell)e^{-\gamma t}(1+t)^m\left[\|g_1\|_{\ell,\infty}+\|g_1\|_{\ell,2}\right].\end{aligned}$$
Applying \eqref{3.7} and \eqref{3.10} gives
$$ \left \|\partial_t^m\partial_x^\ell\mathcal{F}^{-1}\left[e^{\kappa_-(\xi)t}A_{-j}^1(\xi)\varphi_\infty(\xi)\hat{g}_j(\xi)\right]\right\|_2\leqq C(m,\ell,n)t^{-\frac{n}{2}}e^{-ct}\|g_j\|_{(2m+\ell-1-n)^+,2}$$
for $t>0$ and $j=1,2$. By \eqref{3.9} and \eqref{3.10}, we get
\begin{align*}
\left\|\partial_t^m\partial_x^\ell\mathcal{F}^{-1}\left[e^{\kappa_+(\xi)t}A_{+2}^1(\xi)\varphi_\infty(\xi)\hat{g}_j(\xi)\right]\right\|_2\leqq C(m,\ell)e^{-ct}\|g_2\|_{(\ell-1)^+,2}
\end{align*}
for $t>0$, with a suitable constant $c>0$. Next we estimate the $L_\infty$ norm. Put
$$\kappa_-(\xi)=-\delta|\xi|^2+\sigma_-(\xi),$$
First, we consider $\kappa_-$, then
$$\begin{aligned}&\partial_{t}^{m}\partial_{x}^{\alpha}\mathcal{F}^{-1}\left[e^{t\kappa_-(\xi)}A_{-j}^1(\xi)\varphi_\infty(\xi)\hat{g}_j(\xi)\right]\\&=\mathcal{F}^{-1}\left[e^{-\delta|\xi|^2t}\right]*\mathcal{F}^{-1}\left[\kappa_-(\xi)^m\xi^\alpha e^{t\sigma_-(\xi)}A_{-j}^1(\xi)\varphi_\infty(\xi)\hat{g}_j(\xi)\right],\end{aligned}$$
which together with \eqref{3.6} and \eqref{3.7} implies
$$\begin{aligned}&\left\|\partial_t^m\partial_x^\alpha\mathcal{F}\left[e^{t\kappa_-(\xi)}A_{-j}^1(\xi)\varphi_\infty(\xi)\hat{g}_j(\xi)\right]\right\|_\infty\\&\leqq C t^{-\frac{3}{4}}\left\|\mathcal{F}^{-1}\left[\kappa_-(\xi)^m\xi^\alpha e^{t\sigma_-(\xi)}A_{-j}^1(\xi)\varphi_\infty(\xi)\hat{g}_j(\xi)\right]\right\|_2\\&\leqq C(\alpha,m,n)t^{-\left(\frac{3}{4}+\frac{n}{2}\right)}e^{-ct}\|g_j\|_{(2m+|\alpha|-n-1)^+,2}.\end{aligned}$$
By Sobolev's inequality, \eqref{3.9} and \eqref{3.10}, we have
$$\begin{aligned}&\left\|\partial_t^m\partial_x^\ell\mathcal{F}^{-1}\left[e^{t\kappa_+(\xi)}A_{+2}^1(\xi)\varphi_\infty(\xi)\hat{g}_2(\xi)\right]\right\|_\infty\\&\leqq C\left\|\partial_t^m\partial_x^\ell\mathcal{F}^{-1}\left[e^{t\lambda_3(\xi)}A_{+2}^1(\xi)\varphi_\infty(\xi)\hat{g}_2(\xi)\right]\right\|_{2,2}\\&\leqq C(m,\alpha)e^{-ct}\|g_2\|_{\ell+1,2}.\end{aligned}$$
Summing up, we derive
$$\begin{gathered}\begin{aligned}\|\partial_t^m\partial_x^\ell \varphi_\infty n(x,t)\|_2\end{aligned} \leqq C(m,\ell,n)e^{-ct}\left[t^{-\frac{n}{2}}\|G\|_{(2m+\ell-n-1)^{+},2}+\|G\|_{\mathbb{W}_{2}^{\ell,(\ell-1)^{+}}}\right] ;\end{gathered}$$
\begin{align*}
&\|\partial_t^m\partial_x^\ell \varphi_\infty (n(x,t)\|_\infty \leqq C(m,\ell,n)t^{-\left(\frac{n}{2}+\frac{3}{4}\right)}e^{-ct}\|G\|_{(2m+\ell-n-1)^+,2}
+C(m,\ell)e^{-ct}\left[\|G\|_{\mathbb{W}_2^{\ell,\ell+1}}+\|n_0\|_{H^{l+2}}.\right]
\end{align*}
We also write
$$\hat{d}(t,\xi)=\sum_{k={\pm}}\sum_{j=1,2}e^{\kappa_k(\xi)t}A_{kj}^1(\xi)\varphi_\infty(\xi)\hat{g}_j(\xi),$$
and
\begin{align*}
\left|\xi^\alpha\partial_\xi^\alpha A_{k1}^1(\xi)\right|\leq C_\alpha|\xi|^{-1},\ \left|\xi^\alpha\partial_\xi^\alpha A_{+2}^2(\xi)\right|\leq C_\alpha(|\xi|)^2,\ \left|\xi^\alpha\partial_\xi^\alpha A_{-2}^2(\xi)\right|\leq C_\alpha
\end{align*}
for $|\xi|\geqq r_2$ and $\kappa={\pm}$.
Similarly, we have
$$\begin{gathered}
\begin{aligned}\|\partial_t^m\partial_x^\ell \varphi_\infty d\|_2\end{aligned} \leq C(m,\ell,n)e^{-ct}\left[t^{-\frac{n}{2}}\|G\|_{W_2^{(2m+\ell-n-1)^+,(2m+\ell-n)^+}}+\|G\|_{W_2^{(\ell-1)^+,(\ell-2)^+}}\right],  \quad  \forall t>0\end{gathered},$$
\begin{align*}
\|\partial_t^m\partial_x^\ell \varphi_\infty d\|_\infty \leqq C(m,\ell,n)e^{-ct}\left[t^{-\left(\frac{n}{2}+\frac{3}{4}\right)}\|G\|_{\mathbb{W}_2^{(2m+\ell-n-1)^+,(2m+\ell-n)^+}} +\|G\|_{\mathbb{W}_2^{\ell+1,\ell}}\right],\quad\forall t>0.
\end{align*}
Choose $\varphi_M(\xi)\in C_0^\infty(\mathbb{R}^3)$ such that $\varphi_0(\xi)+\varphi_M(\xi)+\varphi_\infty(\xi)=1$ for $\forall\xi\in\mathbb{R}^3$ and Fourier multiplier theorem, we obtain
$$\left\|\partial_t^m\partial_x^\ell\left[\mathcal{F}^{-1}\left[e^{-t\hat{B}(\xi)}\varphi_M(\xi)\hat{G}(\xi)\right]\right]\right\|_2\leqq C(m,\ell)e^{-ct}\|G\|_p$$
for a suitable constant $c>0$.\par
The system \eqref{3.2} is the same as system \eqref{3.1}. Similarly, $D$ be a matrix-valued differential operator given by
$$D:=\begin{bmatrix}0&-\Lambda I_{3\times3}\\\\ \Lambda I_{3\times3}&\mu\Delta I_{3\times3}\end{bmatrix}.$$
Next, analyze the differential operator $D$ in terms of its Fourier expression $\hat{D}(\xi)$,
$$\hat{D}(\xi):=\begin{bmatrix}0&-|\xi| I_{3\times3}\\\\|\xi| I_{3\times3}&-\mu|\xi|^{2} I_{3\times3}\end{bmatrix}.$$
The characteristic polynomial of $\hat{D}(\xi)$ is $\kappa^3(\kappa+\mu|\xi|^2+\frac{|\xi|^2}{\kappa})^3=0$, which implies the eigenvalues are
$$\kappa_{\pm}=-\frac{1}{2}\mu|\xi|^{2}\pm\frac{1}{2}\sqrt{\mu^{2}|\xi|^{4}-4|\xi|^{2}}\quad(triple).$$
The system \eqref{3.2} can also have the same estimates like the system \eqref{3.1}. Combining all the estimates, and
$$\begin{aligned}\|v\|_{L^{2}}\leqq C\|d\|_{L^{2}}+C\|w\|_{L^{2}}.\end{aligned}$$
On the other hand, $\begin{aligned}\|d\|_2=\||\xi|^{-1}\xi^{T}\hat{v}\|_2\leqq C\|\hat{v}\|_2\leqq C\|v\|_2\end{aligned}$. Similarly, variable $\omega$ have the same estimate.
 we have Theorem \ref{th3.1}.
\end{proof}
\section{$L_2$ Decay Estimate in $\Omega$}
In this section we shall prove Theorem \ref{the1.3} by using the cut-off technique based on Lemma\ref{lem2.2} and Theorem \ref{th3.1}. The strategy follows Kobayashi and Shibata \cite{Kobayashi-Shibata1999}. To use Theorem \ref{th3.1}, we construct a suitable extension of $\mathbb{U}_0=(n_0,v_0,E_0)$ to $\mathbb{R}^3$ in the following manner. Let $\mathbb{U}_0^{\prime}$ be the Lions' extension of $\mathbb{U}_0$ such that
$$\mathbb{U}_0^{\prime}|_\Omega=\mathbb{U}_0,\quad\|\mathbb{U}_0^{\prime}\|_{\mathbb{W}_2^{k,\ell}(\mathbb{R}^3)}\leqq C(k,\ell)\|\mathbb{U}_0\|_{\mathbb{W}_2^{k,\ell}}.$$
Put $\mathbb{U}_1(x)=\mathbb{U}_0^\prime=(n_0^\prime,v_0^\prime,E_0^\prime)$.
Then, we have
\begin{align}\label{4.1}
\|\mathbb{U}_1\|_{\mathbb{W}_2^{k,\ell}(\mathbb{R}^3)}\leqq C\|\mathbb{U}_0\|_{\mathbb{W}_2^{k,\ell}}.
\end{align}
Next, we set
$$\mathbb{U}_{c}(t)=(n_c(t),v_c(t),E_c(t))=\mathbb{E}(t)\mathbb{U}_1,\quad\mathrm{in~}(0,\infty)\times\mathbb{R}^3,$$
$$\mathbb{U}(t)=(n(t),v(t),E(t))=e^{t\mathbb{A}}\mathbb{U}_0,\quad\mathrm{~in~}(0,\infty)\times\Omega.$$

We shall start with the following steps.\par
$\emph{Step1.}$ For $t \geqq 1$, we have
\begin{align}\label{1,2,b}
\|\mathbb{U}(t)\|_{1,2,\Omega_b}&\leqq C_bt^{-\frac{3}{2q}}\left(\|\mathbb{U}_0\|_q+\|\mathbb{U}_0\|_{\mathbb{W}_2^{1,0}}\right),\quad\mathbb{U}_0\in\mathbb{W}_2^{1,0}\cap\mathbb{L}_q,
\end{align}
\begin{align}\label{2,2,b}
\|(\mathbb{I}-\mathbb{P})\mathbb{U}(t)\|_{2,2,\Omega_b}&\leqq C_bt^{-\frac{3}{2q}}\left(\|\mathbb{U}_0\|_q+\|\mathbb{U}_0\|_{\mathbb{W}_2^{2,1}}\right),
\quad\mathbb{U}_0\in\mathbb{W}_2^{2,1}\cap\mathbb{L}_q.\end{align}
Put $\begin{aligned}\mathbb{U}_m=(n_m,v_m,E_m)=\mathbb{U}-(1-\varphi)\mathbb{U}_c\end{aligned}$ and then
$$\partial_t\mathbb{U}_m-\mathbb{A}\mathbb{U}_m=\mathbb{F}_m\mathrm{~in~}(0,\infty)\times\Omega,\quad
\mathbb{P}\mathbb{U}_m|_{\partial\Omega}=0,
\quad\mathbb{U}_m(0)=(\varphi n_0,\varphi v_0,\varphi E_0),$$
where $\mathbb{F}_{n}=(f_{m1},f_{m2},f_{m3})$ and
$$\begin{aligned}
f_{m_1}&=(\nabla\varphi)\cdot v_c,\\
(f_{m_2})_i&=\mu[\Delta\varphi v_c^i+2(\nabla_i\varphi:\nabla_i v_c^i)]-(\lambda+\mu)[\nabla_i (\nabla\varphi\cdot v_c)+(\nabla_i \varphi) \mathrm{div}v_c]+(\nabla_i \varphi)n-(E\nabla \varphi)_i,\\
f_{m_3}&=-\nabla\varphi\otimes v_c.\end{aligned}$$
By Duhamel's principle, we obtain
$$\mathbb{U}_m(t)=e^{t\mathbb{A}}\mathbb{U}_m(0)+\int_0^te^{(t-s)\mathbb{A}}\mathbb{F}_m(s)ds.$$
which combined with Lemma \ref{lem2.2} gives
\begin{align*}
\|e^{t\mathrm{A}}\mathbb{U}_m(0)\|_{\mathbb{W}_2^{1,2}(\Omega_b)}\leqq Ct^{-2}\|\mathbb{U}_m(0)\|_{\mathbb{W}_2^{1,0}},
\end{align*}
\begin{align*}
\|e^{t\mathrm{A}}\mathbb{U}_m(0)\|_{\mathbb{W}_p^{2,3}(\Omega_b)}\leqq Ct^{-2}\|\mathbb{U}_m(0)\|_{\mathbb{W}_p^{2,1}}.
\end{align*}
In view of Lemma \ref{lem2.2}, we can write
\begin{align*}
\int_0^te^{(t-s)\mathbb{A}}\mathbb{F}_m(s)ds&\begin{aligned}=\int_0^t\mathbb{T}(t-s)\mathbb{F}_m(s)ds
\end{aligned},
\end{align*}


$$\mathbb{T}(t-s)\mathbb{F}_m(s)=e^{(t-s)\mathbb{A}}(f_{m1}(s),f_{m2}(s),f_{m3}(s)).$$
To estimate $I(t)$, put
$$I(t)=\left\{\int_{t-1}^t+\int_1^{t-1}+\int_0^1\right\}\mathbb{T}(t-s)\mathbb{F}_m(s)ds=\sum_{j=1}^3J_j(t).$$
By Theorem \ref{th3.1} and ({4.1}), we obtain
\begin{align}\label{4.5}
&\|\mathbb{F}_{m}(s)\|_{\mathbb{W}_{p}^{j+1,j}}\nonumber\\
&\leqq C\left[\|(n_c)_s\|_{j+1,p,B_b}+\|\mathbb{E}(s)\mathbb{U}_1\|_{j+1,p,B_b}\right]\nonumber\\
&\leqq Cs^{-\frac{3}{2q}}\left[\|\mathbb{U}_0\|_q+\|\mathbb{U}_0\|_{\mathbb{W}_p^{j+1,j}}\right],\quad s\geqq1,\quad j=0,1.
\end{align}
By \eqref{2.2}, \eqref{2.3} and \eqref{4.5}, we obtain
\begin{align}\label{4.6}
\notag
\|J_1(t)\|_{1,2}&\leqq C\int_{t-1}^t(t-s)^{-\frac{1}{2}}\|(f_{m1}(s),f_{m2}(s),f_{m3}(s))\|_{\mathbb{W}_2^{1,0}}ds\\ \notag
&\leqq C\int_{t-1}^t(t-s)^{-\frac{1}{2}}s^{-\frac{3}{2q}}ds\left[\|\mathbb{U}_0\|_q+\|\mathbb{U}_0\|_{\mathbb{W}_2^{1,0}}\right]\\ \notag
&\leqq Ct^{-\frac{3}{2q}}\left[\left\|\mathbb{U}_0\right\|_q+\left\|\mathbb{U}_0\right\|_{\mathbb{W}_2^{1,0}}\right];\\
\end{align}
\begin{align}\label{4.7}
 \notag
\|(\mathbb{I}-\mathbb{P})J_1(t)\|_{2,2}&\leqq C\int_{t-1}^t(t-s)^{-\frac{1}{2}}\|(f_{m1}(s),f_{m2}(s),f_{m3}(s))\|_{\mathbb{W}_2^{2,1}}ds\\  \notag
&\leqq C\int_{t-1}^t(t-s)^{-\frac{1}{2}}s^{-\frac{3}{2q}}ds\left[\|\mathbb{U}_0\|_q+\|\mathbb{U}_0\|_{\mathbb{W}_2^{2,1}}\right]\\
&\leqq Ct^{-\frac{3}{2q}}\left[\|\mathbb{U}_0\|_q+\|\mathbb{U}_0\|_{\mathbb{W}_2^{2,1}}\right].
\end{align}
By Lemma \ref{lem2.2} and \eqref{4.5}, we obtain
\begin{align*}
 \notag\|J_2(t)\|_{\mathbb{W}_2^{j+1,j+2}(\Omega_b)}&\leqq C\int_1^{t-1}(t-s)^{-2}\|\mathbb{F}_m(s)\|_{\mathbb{W}_2^{j+1,j}}ds\\&\leqq \notag C\int_1^{t-1}(t-s)^{-2}s^{-\frac{3}{2q}}ds\left[\|\mathrm{U}_0\|_q+\|\mathrm{U}_0\|_{\mathbb{W}_2^{j+1,j}}\right]\\&\leqq Ct^{-\frac{3}{2q}}\left[\|\mathbb{U}_0\|_q+\|\mathbb{U}_0\|_{\mathbb{W}_2^{j+1,j}}\right]\quad t\geqq1.
\end{align*}
Since
$$\mathbb{F}_m(s)=(\nabla\varphi\cdot v_c(s),f_{m2}(s),f_{m3}(s)),$$
we have
$$\begin{aligned}J_{3}(t)&=\int_0^1\mathbb{T}(t-s)(\nabla\varphi\cdot v_c(s),f_{m2}(s),f_{m3}(s))ds.\end{aligned}$$
By Theorem \ref{th3.1} and \eqref{4.1}, we obtain
\begin{align*}
\|(\nabla\varphi\cdot v_c(s),f_{m2}(s),f_{m3}(s))\|_{\mathbb{W}_2^{j+1,j}}&\leqq C\|\nabla\varphi \mathbb{E}(s)\mathbb{U}_1\|_{j+1,2,\Omega_b}\\
 &\leqq Cs^{-\frac{1}{2}}\|\mathbb{U}_0\|_{\mathbb{W}_2^{j+1,j}}
\end{align*}
for $j=0,1$ and $0<s\leqq 1$, and thus we get by Lemma \ref{lem2.2} that
\begin{align}\label{4.8}
\|J_3(t)\|_{\mathbb{W}_2^{j+1,j+2}(\Omega_b)}\leqq Ct^{-2}\|\mathbb{U}_0\|_{\mathbb{W}_2^{j+1,j}}\quad j=0,1.
\end{align}
Combining \eqref{4.6}, \eqref{4.7} and \eqref{4.8} implies that
\begin{align*}
\|I(t)\|_{1,2,\Omega_b}&\leqq Ct^{-\frac{3}{2q}}\left[\|\mathbb{U}_0\|_q+\|\mathbb{U}_0\|_{\mathbb{W}_2^{1,0}}\right],\\  \notag
\|(\mathbb{I}-\mathbb{P})I(t)\|_{2,2,\Omega_b}&\leqq Ct^{-\frac{3}{2q}}\left[\|\mathbb{U}_0\|_q+\|\mathbb{U}_0\|_{\mathbb{W}_2^{2,1}}\right]
\end{align*}
for $t\geqq1$.\par

$\emph{Step2.}$ For $t\geqq 1$, we have
\begin{align}\label{4.11}
\|\partial_t\mathbb{U}(t)\|_{1,2,\Omega_b}\leqq Ct^{-\frac{3}{2q}}\left[\|\mathbb{U}_0\|_q+\|\mathbb{U}_0\|_{\mathbb{W}_2^{1,0}}\right]\quad\mathbb{U}_0\in\mathbb{W}_2^{1,0}\cap\mathbb{L}_q,
\end{align}
\begin{align}\label{4.12}
\|\mathbb{P}\mathbb{U}(t)\|_{2,2,\Omega_b}\leqq Ct^{-\frac{3}{2q}}\left[\|\mathbb{U}_0\|_q+\|\mathbb{U}_0\|_{\mathbb{W}_2^{1,0}}\right]\quad\mathbb{U}_0\in\mathbb{W}_2^{1,0}\cap\mathbb{L}_q,
\end{align}
\begin{align}\label{4.13}
\|\mathbb{P}\mathbb{U}(t)\|_{3,2,\Omega_b}\leqq Ct^{-\frac{3}{2q}}\left[\|\mathbb{U}_0\|_q+\|\mathbb{U}_0\|_{\mathbb{W}_2^{2,1}}\right]\quad
\mathbb{U}_0\in\mathbb{W}_2^{2,1}\cap\mathbb{L}_q,\quad\mathbb{P}\mathbb{U}_0|_{\partial\Omega}=0.
\end{align}
By Duhamel's principle and integration by parts, we have
$$\partial_t\mathbb{U}_m(t)=\partial_te^{t\mathbb{A}}\mathbb{U}_m(0)-e^{\frac{1}{2}\mathbb{A}}\mathbb{F}_m(t-\frac{1}{2})
+\int_{t-\frac{1}{2}}^te^{(t-s)\mathbb{A}}\partial_s\mathbb{F}_m(s)ds
+\int_0^{t-\frac{1}{2}}\partial_te^{(t-s)\mathbb{A}}\mathbb{F}_m(s)ds.$$
By Lemma\ref{lem2.2}, Theorem \ref{th3.1} and \eqref{4.1}, one gets
$$\begin{aligned}\|\partial_te^{t\mathbb{A}}\mathbb{U}_m(0)\|_{\mathbb{W}_2^{1,2}(\Omega_b)}&\leqq Ct^{-\frac{3}{2q}}\|\mathbb{U}_0\|_{\mathbb{W}_2^{1,0}},\\\|e^{\frac{1}{2}\mathbb{A}}\mathbb{F}_n(t-\frac{1}{2})\|_{\mathbb{W}_2^{1,2}(\Omega_b)}&\leqq C\|\mathbb{F}_m(t-\frac{1}{2})\|_{\mathbb{W}_2^{1,0}}\leqq C\|\mathbb{E}(t-\frac{1}{2})\mathbb{U}_1\|_{\mathbb{W}_2^{0,1}(\Omega_b)}\\&\leqq Ct^{-\frac{3}{2q}}\left[\left\|\mathbb{U}_0\right\|_q+\left\|\mathbb{U}_0\right\|_{\mathbb{W}_2^{1,0}}\right].\end{aligned}$$
Since
$t-\frac{1}{2}\geqq \frac{1}{2}$, it thus holds that
\begin{align*}
\|\mathbb{E}(t-\frac{1}{2})\mathbb{U}_1\|_{\mathbb{W}_2^{0,1}(\Omega_b)}&\leqq C \|\mathbb{U}_0\|_q+C (t-\frac{1}{2})^{-\frac{3}{2q}}\|\mathbb{U}_0\|_{\mathbb{W}_1^{1,1}}\\
&\leqq t^{-\frac{3}{2q}} \left[\left\|\mathbb{U}_0\right\|_q
+\left\|\mathbb{U}_0\right\|_{\mathbb{W}_2^{1,0}}\right]
\end{align*}
for $\frac{1}{2}\leqq t-\frac{1}{2}\leqq 1$, and
\begin{align*}
\|\mathbb{E}(t-\frac{1}{2})\mathbb{U}_1\|_{\mathbb{W}_2^{0,1}(\Omega_b)}\leqq C(t-\frac{1}{2})^{-\frac{3}{2q}} \left[\left\|\mathbb{U}_0\right\|_q+\left\|\mathbb{U}_0\right\|_{\mathbb{W}_2^{1,0}}\right]
\end{align*}
for $t-\frac{1}{2}\geqq 1$.
Since
\begin{align*}
\|\partial_s\mathbb{F}_m(s)\|_{\mathbb{W}_2^{1,0}}&\leqq C\|\partial_s\mathbb{E}(s)\mathbb{U}_1\|_{\mathbb{W}_2^{0,1}(\Omega_b)}\leqq Cs^{-\left(\frac{3}{2q}+\frac{1}{2}\right)}\left[\|\mathbb{U}_0\|_q+\|\mathbb{U}_0\|_{\mathbb{W}_2^{1,0}}\right]
\end{align*}
for $s\geqq \frac{1}{2}$, by Theorem \ref{th3.1}, \eqref{4.1} and \eqref{2.2}, we obtain
\begin{align*}
\left\|\int_{t-\frac{1}{2}}^te^{(t-s)\mathbb{A}}\partial_s\mathbb{F}_m(s)ds\right\|_{1,2}
&\leqq C\int_{t-\frac{1}{2}}^t(t-s)^{-\frac{1}{2}}s^{-\left(\frac{3}{2q}+\frac{1}{2}\right)}ds\left[\|\mathbb{U}_0\|_q+\|\mathbb{U}_0\|_{\mathbb{W}_2^{1,0}}\right]\\
&\leqq Ct^{-\left(\frac{3}{2q}+\frac{1}{2}\right)}\left[\|\mathbb{U}_0\|_q+\|\mathbb{U}_0\|_{\mathbb{W}_2^{1,0}}\right].
\end{align*}
Since $s \in(0,\frac{1}{2})$, we obtain
\begin{align*}
\|\mathbb{F}_m(s)\|_{\mathbb{W}_2^{1,0}}
&\leqq C\|\mathbb{E}(t)\mathbb{U}_1\|_{\mathbb{W}_2^{0,1}(\Omega_b)}\\
&\leqq C(1+s)^{-\frac{3}{2q}}\left(1+s^{-\frac{1}{2}}\right)\left[\|\mathbb{U}_0\|_q+\|\mathbb{U}_0\|_{\mathbb{W}_2^{1,0}}\right]
\end{align*}
which together with Theorem \ref{th3.1}, \eqref{4.1} and Lemma \ref{lem2.2} leads to
\begin{align*}
\left\|\int_0^{t-\frac{1}{2}}\partial_te^{(t-s)\mathbb{A}}\mathbb{F}_m(s)ds\right\|_{\mathbb{W}_2^{1,2}(\Omega_b)}
\leqq&\ C\int_0^{t-\frac{1}{2}}(t-s)^{-2}(1+s)^{-\frac{3}{2q}}\left(1+s^{-\frac{1}{2}}\right)ds\left[\|\mathbb{U}_0\|_q+\|\mathbb{U}_0\|_{\mathbb{W}_2^{1,0}}\right]\\
\leqq&\ Ct^{-\frac{3}{2q}}\left[\|\mathbb{U}_0\|_q+\|\mathbb{U}_0\|_{\mathbb{W}_2^{1,0}}\right].
\end{align*}
With these estimates at hand, we have
\begin{align}\label{4.14}
\|\partial_t\mathbb{U}_m(t)\|_{1,2,\Omega_b}\leqq Ct^{-\frac{3}{2q}}\left[\|\mathbb{U}_0\|_q+\|\mathbb{U}_0\|_{\mathbb{W}_2^{1,0}}\right].
\end{align}
On the other hand, by Theorem \ref{th3.1} and \eqref{4.1}, we obtain
$$\|\partial_t\mathbb{U}_c(t)\|_{1,2,\Omega_b}\leqq Ct^{-\left(\frac{3}{2q}+\frac{1}{2}\right)}\left[\|\mathbb{U}_0\|_q+\|\mathbb{U}_0\|_{\mathbb{W}_2^{1,0}}\right],$$
which together with \eqref{4.14} implies \eqref{4.11}.\par
By Lemma \ref{lem 7.2} in the Appendix, we have
$$\|\mathbb{PU}(t)\|_{2,2,\Omega_{b-1}}\leqq C\left[\|\partial_t\mathbb{PU}(t)\|_{2,\Omega_b}+\|\mathbb{U}(t)\|_{1,2,\Omega_b}\right].$$
Therefore, by Step 1 and \eqref{4.11}, we have
\begin{align*}
\|\mathbb{PU}(t)\|_{2,2,\Omega_{b-1}}\leqq Ct^{-\frac{3}{2q}}\left[\|\mathbb{U}_0\|_q+\|\mathbb{U}_0\|_{\mathbb{W}_2^{1,0}}\right]\quad t\geqq1.
\end{align*}
Since $b$ is chosen arbitrarily, we get \eqref{4.12}. Applying Proposition Ap.2 with $m=1$, we have
$$\|\mathbb{PU}(t)\|_{3,2,\Omega_{b-1}}\leqq C\left[\|\partial_t\mathbb{PU}(t)\|_{1,2,\Omega_b}+\|\mathbb{U}(t)\|_{2,2,\Omega_b}\right].$$
And then, by Step 1, \eqref{4.11} and \eqref{4.12}, we obtain
$$\|\mathbb{PU}(t)\|_{3,2,\Omega_{b-1}}\leqq Ct^{-\frac{3}{2q}}\left[\|\mathbb{U}_0\|_q+\|\mathbb{U}_0\|_{\mathbb{W}_2^{2,1}}\right]$$
for $\mathbb{U}_0\in\mathbb{W}_2^{2,1}\cap\mathbb{L}_q$. Since $b$ is chosen arbitrarily, we have \eqref{4.13}.\par

$\emph{Step3.}$  We shall estimate $\mathbb{U}(t)$ for $|x|\geqq b$.
Let $\varphi_{\infty}(x)$ be a function in $C^\infty(\mathbb{R}^3)$ such that
$\varphi_\infty(x)=1$

\noindent for $|x|$ $\geqq b$ and $\varphi_\infty(x)=0$ for $|x|\leqq b-1$.
Put $\mathbb{U}_\infty(t)=(n_\infty(t),v_\infty(t),E_\infty(t))=\varphi_\infty\mathbb{U}(t).$ Then,
$$\partial_t\mathbb{U}_\infty-\mathbb{A}\mathbb{U}_\infty=\mathbb{F}_\infty\mathrm{~in~}(0,\infty)\times\mathbb{R}^3\mathrm{~and~}\mathbb{U}_\infty(0)=\varphi_\infty\mathbb{U}_0\mathrm{~in~}\mathbb{R}^3,$$
where $\mathbb{F}_\infty=(f_{\infty\mathrm{l}},f_{\infty2},f_{\infty3})$, and
\begin{align*}
f_{\infty_1}&=\nabla \varphi_\infty\cdot v,\\
(f_{\infty2})_i&=-\mu[\Delta\varphi _\infty v_{c}^i+2\nabla_i\varphi _{\infty}:\nabla_i v^i_c]-(\lambda+\mu)[\nabla_i(\nabla\varphi_{\infty}\cdot v_{c})+\nabla_i\varphi_{\infty}\mathrm{div}v_{c}]+\nabla_i \varphi_\infty n-(E\nabla \varphi_\infty)_i,\\
f_{\infty_3}&=-\nabla \varphi_{\infty}\otimes v.
\end{align*}
By Duhamel's principle, we can write
$$\mathbb{U}_\infty(t)=\mathbb{E}(t)\mathbb{U}_\infty(0)+\mathbb{L}(t),\quad\mathbb{L}(t)=\int_0^t\mathbb{E}(t-s)\mathbb{F}_\infty(s)ds.$$
By \eqref{1,2,b} and Lemma \ref{lem 7.2} in the Appendix below, we have
$$\begin{aligned}\|\mathbb{U}\|_{\infty,\Omega_b}&
\\&\leqq\|\mathbb{U}\|_{1,2,\Omega_b}\\&\leqq Ct^{-\frac{3}{2q}}\left[\|\mathbb{U}_{0}\|_{q}+\|\mathbb{U}_{0}\|_{W_{2}^{1,0}}\right].\end{aligned}$$
By \eqref{2,2,b}, we obtain
$$\begin{aligned}\|\mathbb{U}\|_{W_{\infty}^{0,1}(\Omega_b)}&\leqq\|(\mathbb{I-P)U}\|_{\infty,\Omega_{b}}
+\|\mathbb{PU}\|_{1,\infty,\Omega_{b}}
\\&\leqq Ct^{-\frac{3}{2q}}\left[\|\mathbb{U}_{0}\|_{q}+\|\mathbb{U}_{0}\|_{W_{2}^{1,0}}\right].\end{aligned}$$
By \eqref{2,2,b} and \eqref{4.12}, we obtain
$$\begin{aligned}\|\mathbb{U}\|_{\mathbb{W}_{\infty}^{1.2}(\Omega_b)}&\leqq \|\mathbb{(I-P)U}\|_{1,\infty,\Omega_{b}}+\|\mathbb{PU}\|_{2,\infty,\Omega_{b}}\\&
\leqq Ct^{-\frac{3}{2q}}\left[\|\mathbb{U}_{0}\|_{q}+\|\mathbb{U}_{0}\|_{\mathbb{W}_{2}^{2,1}}\right].\end{aligned}$$

Let $1\leqq q\leqq2$. We estimate the $L_\infty$ norm of $\mathbb{U}(t)$ with $t\geqq1$. Set
$$\sigma=\frac{3}{2}\left(\frac{1}{q}-\frac{1}{2}\right)$$
By Theorem \ref{th3.1}, we have
\begin{align}\label{4.15}
\notag
\|\partial_x^k(\mathbb{I}-\mathbb{P})\mathbb{E}(t)\mathbb{U}_\infty(0)\|_p&\leqq Ct^{-\left(\sigma+\frac{k}{2}\right)}\left[\left\|\mathrm{U}_0\right\|_q+\left\|\mathrm{U}_0\right\|_{\mathbb{W}_p^{k,(k-1)^+}}\right],\\ \notag
\|\partial_x^k\mathbb{PE}(t)\mathbb{U}_\infty(0)\|_p&\leqq Ct^{-\left(\sigma+\frac{k}{2}\right)}\left[\|\mathbb{U}_0\|_q+\|\mathbb{U}_0\|_{\mathbb{W}_p^{(k-1)^+,(k-2)^+}}\right],\\ \notag
\|\mathbb{E}(t)\mathbb{U}_\infty(0)\|_{\mathbb{W}_\infty^{0,1}}&\leqq Ct^{-\frac{3}{2q}}\left[\|\mathbb{U}_0\|_q+\|\mathbb{U}_0\|_{\mathbb{W}_p^{2,1}}\right],\\
\|\mathbb{E}(t)\mathbb{U}_\infty(0)\|_{\mathbb{W}_\infty^{1,2}}&\leqq Ct^{-\frac{3}{2q}}\left[\|\mathbb{U}_0\|_q+\|\mathbb{U}_0\|_{\mathbb{W}_p^{3,2}}\right].
\end{align}
To estimate $\mathbb{L}(t)$, we decompose it as follows:
$$\mathbb{L}(t)=\left\{\int_{t-1}^t+\int_1^{t-1}+\int_0^1\right\}\mathbb{E}(t-s)\mathbb{F}_\infty(s)ds=\sum_{j=1}^3\mathbb{L}_j(t).$$
Since $\|\mathbb{F}_\infty(s)\|_{q}\leqq C\|\mathbb{F}_\infty(s)\|_{2}$ and $\mathbb{F}_\infty$ has the compact support, we deduce from Theorem \ref{th3.1} that
\begin{align*}
\|\partial_x^k\|\mathbb{PL}_1(t)\|_2&\leqq C\int_{t-1}^t(t-s)^{-\frac{1}{2}}\left[\|F_\infty\|_{q}+\|F_\infty\|_{W_2^{(k-1)^+,(k-2)^+}}\right]ds\\&\leqq C\int_{t-1}^t(t-s)^{-\frac{1}{2}}\|F_\infty\|_{(k-1)^+,2}ds,\\
\|\partial_x^{k}(\mathbb{I-P)L}_{1}(t)\|_{2}&\leqq C\int_{t-1}^{t}(t-s)^{-\frac{1}{2}}\left[\|F_{\infty}\|_{q}+\|F_{\infty}\|_{w_{2}^{k,(k-1)^{+}}}\right]ds\\&\leqq C\int_{t-1}^{t}(t-s)^{-\frac{1}{2}}\|F_{\infty}\|_{W_{2}^{k,(k-1)^{+}}}ds,\\
\|\mathbb{L}_{1}(t)\|_{\mathbb{W}_{\infty}^{0,1}}&\leqq C\int_{t-1}^{t}(t-s)^{-\frac{1}{2}}\left[\|F_{\infty}\|_{q}+\|F_{\infty}\|_{W_{2}^{2,1}}\right]ds\\&\leqq C\int_{t-1}^{t}(t-s)^{-\frac{1}{2}}\|F_{\infty}\|_{W_{2}^{2,1}ds},\\
\|\mathbb{L}_{1}\|_{\mathbb{W}_{2}^{1,2}}&\leqq C\int_{t+1}^{t}(t-s)^{-\frac{1}{2}}\left[\|F_{\infty}\|_{q}+\|F_{\infty}\|_{W_{2}^{2,1}}\right]ds\\&\leqq C\int_{t+1}^{t}(t-s)^{-\frac{1}{2}}\|F_{\infty}\|_{W_{2}^{2,3}}ds,
\end{align*}
By Step 1 and Step 2, we obtain
\begin{align}\label{4.16}
\|\mathbb{F}_\infty(s)\|_{\mathbb{W}_2^{j+2,j+1}(\mathbb{R}^3)}\leqq C\|\mathbb{U}(s)\|_{\mathbb{W}_2^{j+1,j+2}(\Omega_b)}\leqq Ct^{-\frac{3}{2q}}\left[\|\mathbb{U}_0\|_q+\|\mathbb{U}_0\|_{\mathbb{W}_2^{j+1,j}}\right]
\end{align}
for $s\geqq1$ and $j=0,1$. Combining these estimate implies that
\begin{align}\label{4.17}
\notag
\|\partial_x^k\mathbb{L}_1(t)\|_2&\leqq Ct^{-\frac{3}{2q}}\left[\|\mathbb{U}_0\|_q+\|\mathbb{U}_0\|_{\mathbb{W}_2^{1,0}}\right],\quad0\leqq k\leqq2,\\ \notag
\|\partial_x^3\mathbb{P}\mathbb{L}_1(t)\|_2&\leqq Ct^{-\frac{3}{2q}}\left[\|\mathbb{U}_0\|_q+\|\mathbb{U}_0\|_{\mathbb{W}_2^{2,1}}\right],\\ \notag
\|\mathbb{L}_1(t)\|_{\mathbb{W}_\infty^{0,1}}&\leqq Ct^{-\frac{3}{2q}}\left[\|\mathbb{U}_0\|_q+\|\mathbb{U}_0\|_{\mathbb{W}_2^{1,0}}\right],\\
\|\mathbb{L}_1(t)\|_{\mathbb{W}_\infty^{1,2}}&\leqq Ct^{-\frac{3}{2q}}\left[\|\mathbb{U}_0\|_q+\|\mathbb{U}_0\|_{\mathbb{W}_2^{2,1}}\right].\end{align}
Set $\ell_{2,q}(t)=t^{-(\sigma+\frac{1}{2})}$, by Theorem \ref{th3.1} and \eqref{4.16}, we obtain
$$\begin{aligned}\|\partial_x^k\mathbb{PL}_2(t)\|_2&\leqq C\int_1^{t-1}(t-s)^{-\frac{3}{2}\left(1-\frac{1}{2}\right)-\frac{k}{2}}\left[\|\mathbb{F}_\infty(s)\|_{1,\mathbb{R}^3}+\|\mathbb{F}_\infty(s)\|_{\mathbb{W}_2^{(k-1)^+,(k-2)^+}(\mathbb{R}^3)}\right]ds\\&\leqq C\int_1^{t-1}(t-s)^{-\frac{3}{4}-\frac{k}{2}}s^{-\frac{3}{2q}}ds\left[\|\mathbb{U}_0\|_q+\|\mathbb{U}_0\|_{\mathbb{W}_2^{1,0}}\right],\quad0\leqq k\leqq3;\end{aligned}$$

$$\begin{aligned}\|\partial_x^k(\mathbb{I}-\mathbb{P})\mathbb{L}_2(t)\|_2&\leqq C\int_1^{t-1}(t-s)^{-\frac{3}{2}\left(1-\frac{1}{2}\right)-\frac{k}{2}}\left[\|\mathbb{F}_\infty(s)\|_{1,\mathbb{R}^3}+\|\mathbb{F}_\infty(s)\|_{\mathbb{W}_2^{k,(k-1)^+}(\mathbb{R}^3)}\right]ds\\&\leqq C\int_1^{t-1}(t-s)^{-\frac{3}{4}-\frac{k}{2}}s^{-\frac{3}{2q}}ds\left[\|\mathbb{U}_0\|_q+\|\mathbb{U}_0\|_{\mathbb{W}_2^{1,0}}\right],\quad0\leqq k\leqq2;\end{aligned}$$

$$\begin{aligned}\|\mathbb{L}_{1}(t)\|_{W_{\infty}^{0,1}(\mathbb{R}^{3})}&\leqq C\int_{t-1}^{t}(t-s)^{-\frac{3}{2}}\left[\|F_{\infty}\|_{1,\mathbb{R}^{3}}+\|F_{\infty}\|_{W_{2}^{2,1},\mathbb{R}^{3}}\right]ds\\&\leq C\int_{t-1}^{t}(t-s)^{-\frac{3}{2}}s^{-\frac{3}{2q}}ds\left[\|\mathbb{U}_0\|_q+\|\mathbb{U}_0\|_{\mathbb{W}^{1,0}_2}\right];\end{aligned}$$

$$\begin{aligned}\|\mathbb{L}_2(t)\|_{\mathbb{W}_\infty^{1,2}(\mathbb{R}^3)}&\leqq C\int_1^{t-1}(t-s)^{-\frac{3}{2}}\left[\|\mathbb{F}_\infty(s)\|_{1,\mathbb{R}^3}+\|\mathbb{F}_\infty(s)\|_{\mathbb{W}_2^{3,2}}\right]ds\\&\leqq C\int_1^{t-1}(t-s)^{-\frac{3}{2}}s^{-\frac{3}{2q}}ds\left[\|\mathbb{U}_0\|_q+\|\mathbb{U}_0\|_{\mathbb{W}_2^{2,1}}\right].\end{aligned}$$

If we divide $\int_1^{t-1}$ by $\int_{\frac{t-1}{2}}^{t-1}$ and $\int_1^{\frac{t-1}{2}}$, we have
$$\begin{aligned}\int_1^{t-1}(t-s)^{-\frac{3}{4}}s^{-\frac{3}{2q}}ds&\leqq Ct^{-\sigma},\\
\int_1^{t-1}(t-s)^{-\frac{3}{4}-\frac{1}{2}}s^{-\frac{3}{2q}}ds&\leqq C\ell_{2q}(t),\\
\int_1^{t-1}(t-s)^{-\frac{3}{4}-\frac{k}{2}}s^{-\frac{3}{2q}}ds&\leqq Ct^{-\frac{3}{2q}},\quad k\geqq2,\\
\int_1^{t-1}(t-s)^{-\frac{3}{2}}s^{-\frac{3}{2q}}ds&\leqq Ct^{-\frac{3}{2q}}.\end{aligned}$$
Therefore, we have
\begin{align}\label{4.18}
\notag
\|\mathbb{L}_2(t)\|_2\leqq Ct^{-\sigma}\left[\|\mathbb{U}_0\|_q+\|\mathbb{U}_0\|_{\mathbb{W}_2^{1,0}}\right],\\ \notag
\|\partial_x\mathbb{L}_2(t)\|_2\leqq C\ell_{2q}(t)\left[\|\mathrm{U}_0\|_q+\|\mathrm{U}_0\|_{\mathbb{W}_2^{1,0}}\right],\\ \notag
\|\partial_x^2\mathbb{L}_2(t)\|_2\leqq Ct^{-\frac{3}{2q}}\left[\|\mathbb{U}_0\|_q+\|\mathbb{U}_0\|_{\mathbb{W}_2^{1,0}}\right],\\ \notag
\|\partial_x^3\mathbb{PL}_2(t)\|_2\leqq Ct^{-\frac{3}{2q}}\left[\|\mathbb{U}_0\|_q+\|\mathbb{U}_0\|_{\mathbb{W}_2^{1,0}}\right],\\ \notag
\|\mathbb{L}_2(t)\|_{\mathbb{W}_\infty^{0,1}}\leqq Ct^{-\frac{3}{2q}}\left[\|\mathbb{U}_0\|_q+\|\mathbb{U}_0\|_{\mathbb{W}_2^{1,0}}\right],\\
\|\mathbb{L}_2(t)\|_{\mathbb{W}_\infty^{1,2}}\leqq Ct^{-\frac{3}{2q}}\left[\|\mathbb{U}_0\|_q+\|\mathbb{U}_0\|_{\mathbb{W}_2^{2,1}}\right].\end{align}
Finally, we estimate $\mathbb{L}_3(t)$. By Theorem \ref{th3.1}, we obtain
$$\begin{aligned}\|\partial_x^k\mathbb{PL}_3(t)\|_2&\leqq C\int_0^1(t-s)^{-\sigma-\frac{k}{2}}\left[\|\mathbb{F}_\infty(s)\|_{q,\mathbb{R}^3}+\|\mathbb{F}_\infty(s)\|_{\mathbb{W}_2^{(k-1)^+,(k-2)^+}(\mathbb{R}^3)}\right]ds\\&\leqq Ct^{-\sigma-\frac{k}{2}}\int_0^1\|\mathrm{U}(s)\|_{\mathbb{W}_2^{(k-2)^+,(k-1)^+}(\Omega_b)}ds;\end{aligned}$$
$$\begin{aligned}\|\partial_x^k(\mathbb{I}-\mathbb{P})\mathbb{L}_3(t)\|_2&\leqq C\int_0^1(t-s)^{-\sigma-\frac{k}{2}}\left[\|\mathbb{F}_\infty(s)\|_{q,\mathbb{R}^3}+\|\mathbb{F}_\infty(s)\|_{\mathbb{W}_2^{k,(k-1)^+}(\mathbb{R}^3)}\right]ds\\&\leqq Ct^{-\sigma-\frac{k}{2}}\int_0^1\|\mathbb{U}(s)\|_{\mathbb{W}_2^{(k-1)^+,k}(\Omega_b)}ds;\end{aligned}$$
$$\begin{aligned}\|\mathbb{L}_3(t)\|_{\mathbb{W}_\infty^{0,1}}&\leqq C\int_0^1(t-s)^{-\frac{3}{2}}\left[\|\mathbb{F}_\infty(s)\|_{1,\mathbb{R}^3}+\|\mathbb{F}_\infty(s)\|_{\mathbb{W}_2^{1,0}(\mathbb{R}^3)}\right]ds\\&\leqq Ct^{-\frac{3}{2}}\int_0^1\|\mathbb{U}(s)\|_{\mathbb{W}_2^{1,2}(\Omega_b)}ds,\end{aligned}$$
$$\begin{aligned}\|\mathbb{L}_3(t)\|_{W_\infty^{1,2}}&\leqq C\int_0^1(t-s)^{-\frac{3}{2}}\left[\|\mathbb{F}_\infty(s)\|_{1,\mathbb{R}^3}+\|\mathbb{F}_\infty(s)\|_{\mathbb{W}_2^{2,1}(\mathbb{R}^3)}\right]ds\\&\leqq Ct^{-\frac{3}{2}}\int_0^1\|\mathbb{U}(s)\|_{\mathbb{W}_2^{2,3}(\Omega_b)}ds.\end{aligned}$$
by \eqref{2.2}, \eqref{2.3}, \eqref{2.4} and \eqref{2.6}, we obtain
$$\begin{aligned}&\|\mathbb{U}(s)\|_{\mathbb{W}_2^{0,1}(\Omega_b)}\leqq Cs^{-\frac{1}{2}}\|\mathbb{U}_0\|_{\mathbb{W}_2^{1,0}},\\&\|\mathbb{U}(s)\|_{\mathbb{W}_2^{1,2}(\Omega_b)}\leqq\|\mathbb{U}_0\|_{\mathbb{W}_2^{1,2}},\quad\mathrm{when~}\mathbb{U}_0\in\mathcal{D}_2(\mathbb{A}),\\&\|\mathbb{U}(s)\|_{\mathbb{W}_2^{2,3}(\Omega_b)}\leqq Cs^{-\frac{1}{2}} \|\mathbb{U}_0\|_{\mathbb{W}_2^{2,1}}.\end{aligned}$$
Combining these estimates implies that
\begin{align}\label{4.19}
\notag
\|\mathbb{L}_3(t)\|_2&\leqq Ct^{-\sigma}\|\mathbb{U}_0\|_{\mathbb{W}_2^{1,0}},\quad\mathbb{U}_0\in\mathbb{W}_p^{1,0},\\ \notag
\|\partial_x\mathbb{L}_3(t)\|_2&\leqq Ct^{-\sigma-\frac{1}{2}}\|\mathbb{U}_0\|_{\mathbb{W}_2^{1,0}},\mathbb{U}_0\in\mathbb{W}_2^{1,0},\\ \notag
\|\partial_x^2\mathbb{PL}_3(t)\|_2&\leqq Ct^{-\sigma-1}\|\mathbb{U}_0\|_{\mathbb{W}_2^{1,0}},\mathbb{U}_0\in\mathbb{W}_2^{1,0},\\ \notag
\|\partial_x^2(\mathbb{I}-\mathbb{P})\mathbb{L}_3(t)\|_2&\leqq Ct^{-\sigma-1}\|\mathbb{U}_0\|_{\mathbb{W}_2^{1,2}},\mathbb{U}_0\in\mathcal{D}_2(\mathbb{A}),\\ \notag
\|\partial_x^3\mathbb{PL}_3(t)\|_2&\leqq Ct^{-\sigma-\frac{3}{2}}\|\mathbb{U}_0\|_{\mathbb{W}_2^{1,2}},\mathbb{U}_0\in\mathcal{D}_2(\mathbb{A}),\\ \notag
\|\mathbb{L}_3(t)\|_{\mathbb{W}_\infty^{0,1}}&\leqq Ct^{-\frac{3}{2}}\|\mathbb{U}_0\|_{\mathbb{W}_2^{1,2}},\quad\mathbb{U}_0\in\mathbb{W}_2^{1,0},\ \mathbb{U}_0\in\mathcal{D}_2(\mathbb{A}),\\
\|\mathbb{L}_3(t)\|_{\mathbb{W}_\infty^{1,2}}&\leqq Ct^{-\frac{3}{2}}\|\mathbb{U}_0\|_{\mathbb{W}_2^{2,1}},\quad\mathbb{U}_0\in\mathcal{D}_2(\mathbb{A}).
\end{align}
Put $D_b=\{x\in\mathbb{R}^3\mid|x|\geqq b\}.$ Since $\mathbb{U}_\infty=\mathbb{U}$ on $D_b$, combining \eqref{4.15}, \eqref{4.17}, \eqref{4.18} and \eqref{4.19} implies that
$$\begin{aligned}\|\mathbb{U}(t)\|_{2,D_b}&\begin{aligned}&\leqq Ct^{-\sigma}\left[\|\mathbb{U}_0\|_q+\|\mathbb{U}_0\|_{\mathbb{W}_2^{1,0}}\right],&\mathbb{U}_0\in\mathbb{W}_2^{1,0}\cap\mathbb{L}_q,\end{aligned}\\\|\partial_x\mathbb{U}(t)\|_{2,D_b}&\leqq C\ell_{2q}(t)\left[\|\mathbb{U}_0\|_q+\|\mathbb{U}_0\|_{\mathbb{W}_2^{1,0}}\right],\mathbb{U}_0\in\mathbb{W}_2^{1,0}\cap\mathbb{L}_q,&\\\|\partial_x^2\mathbb{P}\mathbb{U}(t)\|_{2,D_b}&\leqq Ct^{-\frac{3}{2q}}\left[\|\mathbb{U}_0\|_q+\|\mathbb{U}_0\|_{\mathbb{W}_2^{1,0}}\right],\mathbb{U}_0\in\mathbb{W}_2^{1,0}\cap\mathbb{L}_q,\\\|\partial_x^2(\mathbb{I}-\mathbb{P})\mathbb{U}(t)\|_{2,D_b}&\leqq Ct^{-\frac{3}{2q}}\left[\|\mathbb{U}_0\|_q+\|\mathbb{U}_0\|_{2,2}\right],\quad\mathbb{U}_0\in\mathbb{W}_2^2\cap\mathbb{L}_q\mathrm{~and~}\mathbb{P}\mathbb{U}_0|_{_{\partial\Omega}}=0,\\\|\partial_x^3\mathbb{PU}(t)\|_{2,D_b}&\leqq Ct^{-\frac{3}{2q}}\left[\|\mathbb{U}_0\|_q+\|\mathbb{U}_0\|_{2,2}\right],\quad\mathbb{U}_0\in\mathbb{W}_2^2\cap\mathbb{L}_q\mathrm{~and~}\mathbb{P}\mathbb{U}_0|_{\partial\Omega}=0,\\\|\mathbb{U}(t)\|_{\mathbb{W}_\infty^{0,1}(D_b)}&\leqq Ct^{-\frac{3}{2q}}\left[\|\mathbb{U}_0\|_q+\|\mathbb{U}_0\|_{2,2}\right],\quad\mathbb{U}_0\in\mathbb{W}_2^{2}\cap\mathbb{L}_q\mathrm{~and~}\mathbb{P}\mathbb{U}_0|_{\partial\Omega}=0,\\\|\mathbb{U}(t)\|_{\mathbb{W}_\infty^{1,2}(D_b)}&\leqq Ct^{-\frac{3}{2q}}\left[\|\mathbb{U}_0\|_q+\|\mathbb{U}_0\|_{\mathbb{W}_2^{2,1}}\right],\quad\mathbb{U}_0\in\mathbb{W}_2^{2,1}\cap\mathbb{L}_q\mathrm{~and~}\mathbb{P}\mathbb{U}_0|_{_{\partial\Omega}}=0,\end{aligned}$$
which together with Steps 1 and 2 and Sobolev's imbedding theorem implies the proof of $(A)$ and $(B)$ of Theorem \ref{the1.3}. To prove $(C)$ of Theorem \ref{the1.3}, we think of $\mathbb{U}_(t)$ as $\mathbb{U}(t)=e^{(t-\frac{1}{2})\mathbb{A}}(e^{\frac{1}{2}\mathbb{A}}\mathbb{U}_0)$. By \eqref{2.3} and \eqref{2.4}
$$\|e^{\frac{1}{2}\mathbb{A}}\mathbb{U}_0\|_{\mathbb{W}_2^{2,3}}\leqq C\|\mathbb{U}_0\|_{\mathbb{W}_2^{2,1}}\quad\mathrm{when~}\mathbb{U}_0\in\mathbb{W}_2^{2,1}.$$
Moreover,
$$\left.\mathbb{P}(e^{\frac{1}{2}\mathbb{A}}\mathbb{U}_0)\right|_{\partial\Omega}=0\mathrm{~and~}\|e^{\frac{1}{2}\mathbb{A}}\mathbb{U}_0\|_q\leqq C\|\mathbb{U}_0\|_{\mathbb{W}_q^{1,0}}\mathrm{~when~}\mathbb{U}_0\in\mathbb{W}_q^{1,0}.$$
Therefore, applying $(B)$ to $e^{(t-\frac{1}{2})\mathbb{A}}(e^{\frac{1}{2}\mathbb{A}}\mathbb{U}_0)$, we have $(C)$ of Theorem \ref{the1.3}, which completes the proof of Theorem \ref{the1.3}.
\section{A Proof of Theorem \ref{the1.2}}
In this section, we shall prove Theorem \ref{the1.2}. Consider
\begin{align*}
\begin{cases}
n_t+\mathrm{div}v=f_1,\\
v_t-\mu\Delta v-(\lambda+\mu)\nabla\mathrm{div}v+\nabla n-\mathrm{div}E^{T}=f_2,\\
E_t-\nabla v=f_3,
\end{cases}
\end{align*}
where
\begin{align*}
&f_1=-n\nabla\cdot v-v\cdot\nabla n\\
&(f_2)_i=\begin{aligned}&E_{jk}\nabla_{j}E_{ik}-\frac{n}{1+n}\left(\mu\Delta v_{i}+(\lambda+\mu)\nabla_{i}\mathrm{div}v\right)-v\cdot\nabla v_{i}\\&-\left(\frac{P^{\prime}(n+1)}{(1+n)P^{\prime}(1)}-1\right)\nabla_in.\end{aligned}\\
&f_3=\nabla v E-v \cdot\nabla E
\end{align*}
If we put $\mathbb{U}=(n,v,E)$, $\mathbb{U}_0=(n_0,v_0,E_0)$ and $\mathbb{F}(\mathbb{U})=(f_1,f_2,f_3),$ then the IBVP \eqref{1.1'} is reduced to the following equation:
\begin{align}\label{5.3}
\mathbb{U}_t-\mathbb{A}\mathbb{U}=\mathbb{F}(\mathbb{U}),\quad t>0,\quad\mathbb{P}\mathbb{U}|_{\partial\Omega}=0,\quad\mathbb{U}|_{t=0}=\mathbb{U}_0,
\end{align}
where $\mathbb{F}(\mathbb{U})$ is written symbolically as follows:
$$\begin{aligned}\mathbb{(I-P)F(U)}&=(n\partial_xv,\partial_xn v,v\partial_xE,\partial_xvE)\\\mathbb{PF(U)}&=(a(n)n\partial_x^2\mathbb{PU},b(\mathbb{U})\mathbb{U}\partial_x\mathbb{U},C(n)\partial_xv\partial_xv)\end{aligned}$$
Let $N_k(0,\infty),\ k=1,2$, be the quantity defined in Section 1. We have,
\begin{align}\label{5.4}
\notag
&\|\mathbb{U}(t)\|_{3,2}^2+\|\partial_t\mathbb{U}(t)\|_{\mathbb{W}_2^{2,1}}^2\leqq CN_1(0,\infty)^2,\\
&\|\mathbb{U}(t)\|_{4,2}^2+\|\partial_t\mathbb{U}(t)\|_{\mathbb{W}_2^{3,2}}^2\leqq CN_2(0,\infty)^2.
\end{align}
Next, we state the estimate of nonlinear terms $\mathbb{F}(\mathbb{U})$:
\begin{align}\label{5.8}
\notag
\|\mathbb{F}(\mathbb{U})(s)\|_{\mathbb{W}_2^{1,0}}&\leq C\|\mathbb{U}(s)\|_{2,2}^2,\\ \notag
\|\mathbb{F}(\mathbb{U})(s)\|_{\mathbb{W}_p^{1,0}}&\leqq C\left(\|\mathbb{U}(s)\|_\infty\|\partial_x\mathbb{U}(s)\|_{1,p}+\|\partial_x\mathbb{U}(s)\|_p\|\partial_x\mathbb{P}\mathbb{U}(s)\|_{2,2}\right), \notag
&2<p<\infty,\\\|\mathbb{F}(\mathbb{U})(s)\|_{\mathbb{W}_2^{2,1}}&\leqq C\left(\|\mathbb{U}(s)\|_\infty+\|\partial_x\mathbb{U}(s)\|_4\right)\|\partial_x\mathbb{U}(s)\|_{2,2},\\ \notag
\|\partial_s\mathbb{F}(\mathbb{U})(s)\|_{\mathbb{W}_2^{1,0}}&\leqq C\left(\|\partial_s\mathbb{U}(s)\|_{2,2}\|\mathbb{U}(s)\|_\infty+\|\partial_s\mathbb{U}(s)\|_{1,2}\|\partial_x\mathbb{U}(s)\|_{2,2}\right).
\end{align}
In fact, to get some estimate of nonlinear terms $\mathbb{F}(\mathbb{U})$, we use H\"{o}lder's inequality and the inequalities:
\begin{align}\label{5.9}
\notag
&\|u\|_p\leqq C\|u\|_{1,2}\quad2\leqq p\leqq6,\\
&\|u\|_p\leqq C\|u\|_{2,2}\quad6<p\leqq\infty.
\end{align}

By choosing $\|(n_0,v_0,E_0)\|_{k+2,2}=\|(\rho_0-1,u_0,F_0-I)\|_{k+2,2}$ small enough, we can make $N_k(0,+\infty)$ as small as we want. Without loss of generality we may assume that $N_k(0,+\infty)\leqq 1,k=1,2$. \par
To prove Theorem \ref{the1.2}, we reduce the problem \eqref{5.3} to the integral equation:
$$\mathbb{U}(t)=e^{t\mathbb{A}}\mathbb{U}_0+\mathbb{N}(t),\quad\mathbb{N}(t)=\int_0^te^{(t-s)\mathbb{A}}\mathbb{F}(\mathbb{U})(s)ds.$$
$\emph{Step1.}$ Put
$$M_2(t)=\sup_{0\leqq s\leqq t}(1+s)^{\frac{3}{4}}\left(\|\mathbb{U}(s)\|_{\mathbb{W}_2^{2,3}}+\|\partial_s\mathbb{U}(s)\|_{1,2}\right).$$
Then, there exists an $\epsilon_{2}>0$ such that if $N_2(0,\infty)\leqq\epsilon_2$, then
\begin{align}\label{5.5}
M_2(t)\leqq C\left(N_1(0,\infty)+\|\mathbb{U}_0\|_1+\|\mathbb{U}_0\|_{2,2}\right).
\end{align}
First, we consider $0\leqq t \leqq 2$, clear,
\begin{align}\label{5.6}
M_2(t)\leqq CN_1(0,\infty).
\end{align}
Then, we consider $t\leqq 2$. By Theorem \ref{the1.3} $(A)$ and $(B)$, we obtain
$$\|e^{t\mathrm{A}}\mathbb{U}_0\|_2\leqq Ct^{-\frac{3}{4}}\left(\|\mathbb{U}_0\|_1+\|\mathbb{U}_0\|_{\mathbb{W}_2^{1,0}}\right),$$
$$\|\partial_xe^{t\mathbb{A}}\mathbb{U}_0\|_2\leqq Ct^{-\frac{5}{4}}\left(\left\|\mathbb{U}_0\right\|_1+\left\|\mathbb{U}_0\right\|_{\mathbb{W}_2^{1,0}}\right),$$
$$\|\partial_x^2\mathbb{P}e^{t\mathbb{A}}\mathbb{U}_0\|_2\leqq Ct^{-\frac{3}{2}}\left(\|\mathbb{U}_0\|_1+\|\mathbb{U}_0\|_{\mathbb{W}_2^{1,0}}\right),$$
\begin{align}\label{5.7}
\|\partial_x^2(\mathbb{I}-\mathbb{P})e^{t\mathbb{A}}\mathbb{U}_0\|_2+\|\partial_x^3\mathbb{P}e^{t\mathbb{A}}\mathbb{U}_0\|_2+\|\partial_x\partial_te^{-t\mathbb{A}}\mathbb{U}_0\|_2
\end{align}
\begin{align*}
\leqq Ct^{-\frac{3}{2}}\left(\|\mathbb{U}_0\|_1+\|\mathbb{U}_0\|_{2,2}\right)
\end{align*}
for $\left.\mathbb{PU}_0\right|_{\partial\Omega}=0.$ We divide $\mathbb{N}(t)$ into two part as follow:
$$\mathbb{N}(t)=\left\{\int_{t-1}^t+\int_0^{t-1}\right\}e^{(t-s)\mathbb{A}}\mathbb{F}(\mathbb{U})(s)ds=I(t)+II(t).$$
Next, we consider $II(t)$
$$II(t)=\int_{0}^{t-1}e^{(t-s)\mathbb{A}}\mathbb{F(U)}(s)ds$$
and
$$\partial_ tII(t)=\int_{0}^{t-1}\partial_ te^{(t-s)\mathbb{A}}\mathbb{F(U)}(s)ds+e^{\mathbb{A}}\mathbb{F(U)}(t-1).$$
By Theorem \ref{the1.3} $(A)$ and $(C)$ with $q = 2$, we obtain
$$\begin{aligned}&\|II(t)\|_{\mathbb{W}_2^{2,3}}+\|\partial_tII(t)\|_{1,2}\\&\leqq C\left\{\|e^{-\mathbb{A}}\mathbb{F}(\mathbb{U})(t-1)\|_{1,2}+\int_0^{t-1}(t-s)^{-\frac{3}{4}}\left[\|\mathbb{F}(\mathbb{U})(s)\|_1+\|\mathbb{F}(\mathbb{U})(s)\|_{\mathbb{W}_2^{2,1}}\right]ds\right\}.\end{aligned}$$
By H\"{o}lder's inequality, \eqref{5.8} and \eqref{5.9}, we obtain
\begin{align}\label{5.10}
\notag
\|\mathbb{F(U)}\|_1&\leqq\|n\partial_ xv\|_1+\|\partial_ xnv\|_1+\|v\partial _xE\|_1+\|\partial _xvE\|_1\\ \notag
&+\|n\partial_ x^2\mathbb{PU}\|_1+\|\mathbb{U}\partial_ x\mathbb{U}\|_1+\|\partial_ xv\partial _xv\|_1\\ \notag
&\leqq\|\mathbb{U}\|_{1,2}\|\partial _x\mathbb{U}\|_{1,2}\\
\|\mathbb{F}(\mathbb{U})(s)\|_{\mathbb{W}_2^{2,1}}&\leqq C\|\mathbb{U}(s)\|_{2,2}\|\partial_x\mathbb{U}(s)\|_{2,2}.
\end{align}
By \eqref{2.2} and \eqref{5.8}, we obtain
$$\|e^{\mathbb{A}}\mathbb{F}(\mathbb{U})(t-1)\|_{1,2}\leqq C\|\mathbb{F}(\mathbb{U})(t-1)\|_{\mathbb{W}_{2}^{1,0}}\leqq C\|\mathbb{U}(t-1)\|_{2,2}^2.$$
Combining these estimations, we have
\begin{align}\label{5.11}
\notag
&\|II(t)\|_{\mathbb{W}_2^{2,3}}+\|\partial_tII(t)\|_{1,2}\\ \notag
&\leqq C\left[\int_0^{t-1}(t-s)^{-\frac{3}{4}}(1+s)^{-\frac{3}{4}}\|\partial_x\mathbb{U}(s)\|_{2,2}dsM_2(t)+\|\mathbb{U}(t-1)\|_{2,2}^2\right]\\ \notag
&\leqq C\left(\int_0^{t-1}(t-s)^{-\frac{3}{2}}(1+s)^{-\frac{3}{2}}ds\right)^{\frac{1}{2}}\left(\int_0^{t-1}\|\partial_x\mathbb{U}(s)\|_{2,2}^2ds\right)^{\frac{1}{2}}M_2(t).\\ \notag
&+C(1+t)^{-\frac{3}{4}}M_2(t)N_1(0,\infty)\\
&\leqq C(1+t)^{-\frac{3}{4}}N_1(0,\infty)M_2(t).
\end{align}
On the other hand, by \eqref{2.2} and \eqref{5.8}, we obtain
$$\begin{aligned}\|I(t)\|_{1,2}&\leqq C\int_{t-1}^t(t-s)^{-\frac{1}{2}}\|\mathbb{F}(\mathbb{U})(s)\|_{\mathbb{W}_2^{1,0}}ds\\&\leqq C\int_{t-1}^t(t-s)^{-\frac{1}{2}}(1+s)^{-\frac{3}{4}}dsN_1(0,\infty)M_2(t)\\&\leqq C(1+t)^{-\frac{3}{4}}N_1(0,\infty)M_2(t),\end{aligned}$$
which together with \eqref{5.11}, \eqref{5.6} and \eqref{5.7} implies that
\begin{align}\label{5.12}
\|\mathbb{U}(t)\|_{1,2}\leqq C(1+t)^{-\frac{3}{4}}\left[N_1(0,\infty)M_2(t)+N_1(0,\infty)+\|\mathbb{U}_0\|_1+\|\mathbb{U}_0\|_{2,2}\right].
\end{align}
By integration by parts, we obtain
\begin{align}\label{5.13}
\partial_tI(t)=\int_{t-1}^te^{(t-s)\mathbb{A}}\partial_s\mathbb{F}(\mathbb{U})(s)ds,
\end{align}
and therefore by \eqref{2.2}, we obtain
$$\|\partial_tI(t)\|_{1,2}\leqq C\int_{t-1}^t(t-s)^{-\frac{1}{2}}\|\partial_s\mathbb{F}(\mathbb{U})(s)\|_{\mathbb{W}_2^{1,0}}ds.$$
By \eqref{5.8}, \eqref{5.9} and \eqref{5.4}, we obtain
$$\begin{aligned}\|\partial_s\mathbb{F}(\mathbb{U})(s)\|_{\mathbb{W}_2^{1,0}}&\leqq C\left[\|\partial_s\mathbb{U}(s)\|_{2,2}\|\mathbb{U}(s)\|_{2,2}+\|\partial_s\mathbb{U}(s)\|_{1,2}\|\mathbb{U}(s)\|_{3,2}\right]\\&\leqq CN_2(0,\infty)(1+s)^{-\frac{3}{4}}M_2(t),\quad0\leqq s\leqq t.\end{aligned}$$
Combining these estimations, we have
$$\|\partial_tI(t)\|_{1,2}\leqq C(1+t)^{-\frac{3}{4}}N_2(0,\infty)M_2(t),$$
which together with \eqref{5.11}, \eqref{5.6} and \eqref{5.7} implies that
\begin{align}\label{5.14}
\|\partial_t\mathbb{U}(t)\|_{1,2}\leqq C(1+t)^{-\frac{3}{4}}\left[N_2(0,\infty)M_2(t)+N_1(0,\infty)+\|\mathbb{U}_0\|_1+\|\mathbb{U}_0\|_{2,2}\right].
\end{align}
By \eqref{2.2}, \eqref{2.3}, \eqref{5.10} and \eqref{5.4}, we obtain
$$\begin{aligned}\|(\mathbb{I}-\mathbb{P})I(t)\|_{2,2}&\leqq C\int_{t-1}^t(t-s)^{-\frac{1}{2}}\|\mathbb{F}(\mathbb{U})(s)\|_{\mathbb{W}_2^{2,1}}ds\\&\leqq CN_1(0,\infty)\int_{t-1}^t(t-s)^{-\frac{1}{2}}(1+s)^{-\frac{3}{4}}dsM_2(t)\\&\leqq C(1+t)^{-\frac{3}{4}}N_1(0,\infty)M_2(t),\end{aligned}$$
which together with \eqref{5.11}, \eqref{5.6} and \eqref{5.7} implies that
\begin{align}\label{5.15}
\|(\mathbb{I}-\mathbb{P})\mathbb{U}(t)\|_{2,2}\leqq C(1+t)^{-\frac{3}{4}}\left[N_1(0,\infty)M_2(t)+N_1(0,\infty)+\|\mathbb{U}_0\|_1+\|\mathbb{U}_0\|_{2,2}\right].
\end{align}
Applying Lemma\ref{lem 7.3} in the Appendix to the second equation, we have
$$\|v\|_{2+k,2}\leqq C\left[\|\mathbb{U}(t)\|_{1+k,2}+\|\partial_t\mathbb{U}(t)\|_{k,2}+\|\mathbb{F}(\mathbb{U})(t)\|_{k,2}\right],\quad k=0,1.$$
For $k=1$, we use $\|v\|_{2,2}$ and $\|\mathbb{(I-P)(U)}\|_{2,2}$ .\par
Together with \eqref{5.15}, \eqref{5.14} and \eqref{5.12} implies that
$$\|\mathbb{PU}(t)\|_{3,2}\leqq C(1+t)^{-\frac{3}{4}}\left[N_2(0,\infty)M_2(t)+N_1(0,\infty)+\|\mathbb{U}\|_1+\|\mathbb{U}_0\|_{2,2}\right],$$
Therefore, we have
$$M_2(t)\leqq C\left[N_2(0,\infty)M_2(t)+N_1(0,\infty)+\|\mathbb{U}_0\|_1+\|\mathbb{U}_0\|_{2,2}\right],$$
which implies \eqref{5.5} provided that $CN_2(0,\infty)<1$.\par
$\emph{Step2.}$ Put
$$M_\infty(t)=\sup_{0\leqq s\leqq t}(1+s)^{\frac{3}{2}}\|\mathrm{U}(s)\|_\infty.$$
Then, there exists an $\epsilon_{\infty}>0$ such that if $N_1(0,\infty)\leqq\epsilon_{\infty}$ then
\begin{align}\label{5.16}
M_\infty(t)\leqq C[N_2(0,\infty)^2+M_2(t)N_1(0,\infty)+M_2(t)^2+N_1(0,\infty)+\|\mathbb{U}_0\|_1+\|\mathbb{U}_0\|_{2,2}].
\end{align}
When $0\leqq t\leqq2$, by Sobolev's inequality and \eqref{5.4}, we obtain
\begin{align}\label{5.17}
M_{\infty}(t)\leqq C\|\mathbb{U}\|_{\infty}\leqq\|\nabla^{2}\mathbb{U}\|_{2}^{\frac{3}{4}}\|\mathbb{U}\|_{2}^{\frac{1}{4}}\leqq\|\mathbb{U}\|_{2,2}
\end{align}
and therefore we consider the case when $t\geqq2$ below. By Theorem \ref{the1.3}$(A)$ with $q=1$
\begin{align}\label{5.18}
\|e^{-t\mathbb{A}}\mathbb{U}_0\|_{\mathbb{W}_\infty^{0,1}}\leqq Ct^{-\frac{3}{2}}\left[\|\mathbb{U}_0\|_1+\|\mathbb{U}_0\|_{2,2}\right].
\end{align}
By Sobolev's inequality and \eqref{2.2}, we obtain
$$\begin{aligned}
\|I(t)\|_\infty&\leqq C\|I(t)\|_{1,4}\\&\leqq C\int_{t-1}^t(t-s)^{-\frac{1}{2}}\|\mathbb{F}(\mathbb{U})(s)\|_{\mathbb{W}_4^{1,0}}ds.
\end{aligned}$$
For $\|\mathbb{F(U)}(s)\|_{W_{4}^{1,0}}$
\begin{align}\label{5.19}
\|\mathbb{F(U)}(s)\|_{W_{4}^{1,0}}&\leqq C[\|\mathbb{U}(s)\|_{\infty}\|\partial_{x}\mathbb{U}(s)\|_{1,4}+\|\partial_{x}\mathbb{U}(s)\|_{4}\|\partial_{x}\mathbb{PU}(s)\|_{2,2}]\\ \notag
&\leqq C[\|\mathbb{U}(s)\|_{\infty}\|\mathbb{U}(s)\|_{3,2}+\|\mathbb{U}(s)\|^2_{\mathbb{W}_{2}^{2,3}}].
\end{align}
Combining these estimations and \eqref{5.4}, we have
\begin{align}\label{5.20}
\|I(t)\|_\infty\leqq C(1+t)^{-\frac{3}{2}}\left[N_1(0,\infty)M_\infty(t)+M_2(t)^2\right].
\end{align}
By Theorem \ref{the1.3}$(A)$ with $q=1$, \eqref{5.19} and \eqref{5.10}, we obtain
$$\|II(t)\|_{\infty}\leqq C\int_{0}^{t-1}(t-s)^{-\frac{3}{2}}[\|\mathbb{F(U)}(s)\|_{1}+\mathbb{F(U)}(s)\|_{2,2}]ds$$
and
$$\|\mathbb{F(U)}(s)\|_{2,2}\leqq M_{\infty}(t)N_{2}(0,\infty)+M_{2}(t)N_{1}(0,\infty)+N_{2}(0,\infty)^{2},$$
which together with \eqref{5.20}, \eqref{5.18} and \eqref{5.17} implies that
$$M_\infty(t)\leqq C[N_2(0,\infty)M\infty(t)(t)+N_2(0,\infty)^2+M_2(t)N_1(0,\infty)+M_2(t)^2+N_1(0,\infty)+\|\mathbb{U}_0\|_1+\|\mathbb{U}_0\|_{2,2}].$$
If $CN_2(0,\infty) < 1$, then we have \eqref{5.16}.\par
$\emph{Step3.}$ Put
$$\mathcal{M}_2(t)=\sup_{0\leqq s\leqq t}(1+s)^{\frac{5}{4}}\left[\|\partial_s\mathbb{U}(s)\|_{1,2}+\|\partial_x\mathbb{U}(s)\|_{\mathbb{W}_2^{1,2}}\right].$$
Then, there exists an $\epsilon_2^{\prime}>0$ such that if $N_2(0,\infty)\leqq\epsilon_2^{\prime}$, then
\begin{align}\label{5.23}
\notag
\mathcal{M}_2(t)&\leqq C\Big[M_\infty(t)+M_\infty(t)N_1(0,\infty)+M_2(t)^2\\&+N_1(0,\infty)+N_1(0,\infty)^2+\|\mathbb{U}_0\|_1+\|\mathbb{U}_0\|_{2,2}\Big].
\end{align}
In view of \eqref{5.6}, we may concentrate on the case when $t \geqq 2$ only. By \eqref{5.13}, \eqref{2.2}, \eqref{5.4} and \eqref{5.8},
\begin{align*}
\|\partial_tI(t)\|_{1,2}\leqq C(1+t)^{-\frac{5}{4}}\Big[N_2(0,\infty)M_\infty(t)+N_1(0,\infty)\mathcal{M}_2(t)\Big].
\end{align*}
By \eqref{2.2}, \eqref{5.8} and Step 1, we obtain
$$\|e^{\mathrm{A}}\mathbb{F}(\mathbb{U})(t-1)\|_{1,2}\leqq C\|\mathbb{F}(\mathbb{U})(t-1)\|_{\mathbb{W}_2^{1,0}}\leqq C(1+t)^{-\frac{3}{2}}M_2(t)^2.$$
By H\"{o}lder's inequality with the exponent: $\frac{5}{6}=\frac{1}{2}+\frac{1}{3}$ and \eqref{5.9}, we obtain
\begin{align}\label{5.25}
\|\mathbb{F(U)}\|_{\mathbb{W}^{1,0}_{6/5}}\leqq\|\mathbb{U}(s)\|_{2,2}^2.
\end{align}
By \eqref{5.4} and \eqref{5.8}, we obtain
$$\begin{aligned}
&\|\partial_t II(t)\|_{1,2}\\
&\leqq C\left[\|e^{-\mathbb{A}_1}\mathbb{F}(\mathbb{U})(t-1)\|_{1,2}+\int_0^{t-1}(t-s)^{-\frac54}\left(\|\mathbb{F}(\mathbb{U})(s)\|_{\mathbb{W}_{6/5}^{1,0}}+\|\mathbb{F}(\mathbb{U})(s)\|_{\mathbb{W}_2^{2,1}}+\|\mathbb{F}(\mathbb{U})(s)\|_1\right)ds\right].\\
&\leqq C\|e^{-\mathbb{A}}\mathbb{F(U)}(t-1)\|_{1,2}+\int_{0}^{t-1}(t-s)^{-\frac{5}{4}}(M_{2}(t)^{2}+M_{\infty}(t)N_{1}(0,\infty)+N_{1}(0,\infty)^{2})ds.
\end{aligned}
$$
By \eqref{2.2}, \eqref{5.8}, we obtain
$$\|(\mathbb{I}-\mathbb{P})I(t)\|_{2,2}\leqq C(1+t)^{-\frac{3}{2}}N_1(0,\infty)\left(M_\infty(t)+N_1(0,\infty)\right).$$
By Theorem \ref{the1.3} $(C)$ with $q =\frac{6}{5}$, \eqref{5.25} and \eqref{5.8}, we obtain
$$\|\partial_x^2(\mathbb{I}-\mathbb{P})II(t)\|_2\leqq C(1+t)^{-\frac{5}{4}}\left[N_1(0,\infty)\left(M_\infty(t)+N_1(0,\infty)\right)+M_2(t)^2\right].$$
Combining these estimations, \eqref{5.6} and \eqref{5.7}, we have
\begin{align}\label{5.26}
\notag
&\|\partial_t\mathbb{U}(t)\|_{1,2}+\|\partial_x^2(\mathbb{I}-\mathbb{P})\mathbb{U}(t)\|_2\\ \notag
&\leqq C(1+t)^{-\frac{5}{4}}[N_1(0,\infty)\mathcal{M}_2(t)+N_2(0,\infty)M_\infty(t)+M_\infty(t)\\
&+M_2(t)^2+N_1(0,\infty)+N_1(0,\infty)^2+\left\|\mathbb{U}_0\|_1+\|\mathbb{U}_0\|_{2,2}\right],\quad\forall t>0,
\end{align}
where we have used the fact that $N_k(0,\infty)\leqq1,k=1,2.$  By Proposition in the Appendix, we have
\begin{align}
\notag
\|\partial_x^{k+2}\mathbb{P}\mathbb{U}(t)\|_2&\leqq C\Big[\|\partial_t\mathbb{U}(t)\|_{k,2}+\|\mathbb{P}\mathbb{F}(\mathbb{U})(t)\|_{k,2}+\|\partial_x(\mathbb{I}-\mathbb{P})\mathbb{U}(t)\|_{k,2}\Big]\\
&+\|\partial_x\mathbb{P}\mathbb{U}(t)\|_{k,2}+\|\mathbb{P}\mathbb{U}(t)\|_{2,\Omega_b} \Big] ,\quad k=0,1.
\end{align}
Since $\mathbb{PU}(t)|_{\partial\Omega}=0,$ by Poincar\'{e}'s inequality, we obtain
\begin{align}
\|\mathbb{P}\mathbb{U}(t)\|_{2,\Omega_b}\leqq C\|\partial_x\mathbb{P}\mathbb{U}(t)\|_2\leqq C(1+t)^{-\frac{5}{4}}\left[\|\mathbb{U}_0\|_1+\|\mathbb{U}_0\|_{\mathbb{W}_2^{1,0}}+M_2(t)^2\right].
\end{align}
By \eqref{5.8}, \eqref{5.4}, Step 2, we obtain
\begin{align}\label{5.29}
\|\mathbb{PF}(\mathbb{U})(t)\|_{1,2}\leqq C(1+t)^{-\frac{3}{2}}N_1(0,\infty)\left(M_\infty(t)+N_1(0,\infty)\right).
\end{align}
Combining \eqref{5.26}-\eqref{5.29}, we obtain
\begin{align}\label{5.30}
\notag
\mathcal{M}_2(t)&\leqq C\Big[N_1(0,\infty)\mathcal{M}_2(t)+M_\infty(t)+M_\infty(t)N_1(0,\infty)+M_2(t)^2\\&+N_1(0,\infty)+N_1(0,\infty)^2+\|\mathbb{U}_0\|_1+\|\mathbb{U}_0\|_{2,2}\Big].
\end{align}
If $CN_1(0,\infty)<1$, then \eqref{5.23} follows from \eqref{5.30}. Combining all the steps, we have
proved Theorem \ref{the1.2}, which completes the proof of Theorem \ref{the1.2}.
\section{Appendix. A Priori Estimate of an Elliptic Operator}
By Agmon-Douglis-Nirenberg \cite{Agmon-Douglis-Nirenberg1959-1,Agmon-Douglis-Nirenberg1959-2} , we know the following estimate.
\begin{Lemma}\label{lem 7.1}
Let D be a bounded domain with smooth boundary $\partial D$. Let $1<p<\infty$ and let $m$ be an integer such that $m\geqq0$. Suppose that $u\in W_p^{m+2}(D)$ and $f\in  W_p^m(\Omega)$ satisfy the equation:
\begin{align*}
-\alpha\Delta u-\beta\nabla\operatorname{div}u=f\operatorname{in}D \quad\left.and\quad u \right|_{\partial D}=0,
\end{align*}
where $\alpha>0$ and $\alpha+\beta>0$. Then, the following estimate holds:
$$\|u\|_{m+2,p,D}\leqq C(m,p,b)\left[\|f\|_{m,p,D}+\|u\|_{p,D}\right].$$
By using the cut-off function $\varphi$ we can deduce the following proposition from Proposition1 immediately.
\end{Lemma}
\begin{Lemma}\label{lem 7.2}
Let $b$ be an arbitrary number such that $B_{b-3}\supset\mathcal{O}.$ Let $1<p<\infty$ and let $m$ be an integer such that $m\geqq0$. Suppose that $u\in W_p^{m+2}(\Omega)$ and $f\in  W_p^m(\Omega)$ satisfy the equation:
\begin{align*}
-\alpha\Delta u-\beta\nabla\operatorname{div}u=f\, \text{in}\,\Omega\quad \text{and}\quad u|_{\partial\Omega}=0,
\end{align*}
where $\alpha>0$ and $\alpha+\beta>0$. Then, the following estimate holds:
$$\|u\|_{m+2,p,\Omega_{b-1}}\leqq C(m,p,b)\left[\|f\|_{m,p,\Omega_b}+\|u\|_{p,\Omega_b}\right].$$
By Fourier transform, we can reduce the formula: $-\alpha\Delta u-\beta\nabla\operatorname{div}u=f$ in $\mathbb{R}^3$ to the formula:
$$\left(\alpha|\xi|^2\delta_{ij}+\beta\xi_i\xi_j\right)\hat{u}(\xi)=\hat{f}(\xi)\quad\mathrm{in~}\mathbb{R}^3,$$
Since
$$\det\left(\alpha|\xi|^2\delta_{ij}+\beta\xi_i\xi_j\right)=\alpha^2(\alpha+\beta)|\xi|^6,$$
by Fourier multiplier theorem we see that
\begin{align}\label{ap3}
\left\|\partial_x^{m+2}\mathcal{F}^{-1}\left[\left(\alpha|\xi|^2\delta_{ij}+\beta\xi_i\xi_j\right)^{-1}\hat{f}(\xi)\right]\right\|_{p,\mathbb{R}^3}\leqq C(m,p)\|f\|_{m,p,\mathbb{R}^3}.
\end{align}
By using the cut-off function, we can deduce the following proposition from \eqref{ap3} and Lemma \ref{lem 7.2}.
\end{Lemma}
\begin{Lemma}\label{lem 7.3}
Let $1<p<\infty$ and let $m$ be an integer such that $m\geqq0$. Suppose that $u\in W_p^{m+2}(\Omega)$ and $f\in  W_p^m(\Omega)$ satisfy the equation:
\begin{align*}
-\alpha\Delta u-\beta\nabla\operatorname{div}u=f\, \text{in}\,\Omega\quad \text{and}\quad u|_{\partial\Omega}=0,
\end{align*}
where $\alpha>0$ and $\alpha+\beta>0$. Then, the following estimate holds:
$$\|\partial_x^{m+2}u\|_p\leqq C(m,p)\left\{\|f\|_{m,p}+\|u\|_{1,p,\Omega_b}\right\},$$
where b is the same as in Lemma \ref{lem 7.2}.
\end{Lemma}

\section*{Acknowledgements}
This work was partially supported by National Key R\&D Program of China (No. 2021YFA1002900) and Guangzhou City Basic and Applied Basic Research Fund (No. 2024A04J6336).

\bigskip

{\bf Data Availability:} Data sharing is not applicable to this article.

\bigskip

{\bf Conflict of Interest:} The authors declare that they have no conflict of interest.


\end{document}